\documentstyle[amscd,amssymb,verbatim,12pt]{amsart}
\newcommand{\BU}{\mbox{\rm BU}}
\newcommand{\C}{{\Bbb C}}

\newcommand{\ch}{\mbox{\rm ch}}

\newcommand{\CS}{\mbox{\rm CS}}

\newcommand{\End}{\mbox{\rm End}}

\newcommand{\HH}{\mbox{\rm H}}

\newcommand{\Image}{\mbox{\rm Im}}

\newcommand{\Ind}{\mbox{\rm Ind}}

\newcommand{\Isom}{\mbox{\rm Isom}}
\newcommand{\KK}{\mbox{\rm K}}
\newcommand{\Ker}{\mbox{\rm Ker}}

\newcommand{\LIM}{\mbox{\rm LIM}}

\newcommand{\pt}{\mbox{\rm pt}}
\newcommand{\Q}{{\Bbb Q}}
\newcommand{\R}{{\Bbb R}}

\newcommand{\rk}{\mbox{\rm rk}}

\newcommand{\Tr}{\mbox{\rm Tr}}

\newcommand{\U}{\mbox{\rm U}}

\newcommand{\Z}{{\Bbb Z}}
\theoremstyle{plain}

\newtheorem{lemma}{Lemma}
\newtheorem{theorem}{Theorem}
\newtheorem{proposition}{Proposition}

\numberwithin{equation}{section}
\renewcommand{\rm}{\normalshape}
\begin{document}
\title{Higher-Degree Analogs of the Determinant Line Bundle}
\author{John Lott}
\address{Department of Mathematics\\
University of Michigan\\
Ann Arbor, MI  48109-1109\\
USA}
\email{lott@@math.lsa.umich.edu}
\thanks{Research supported by NSF grant DMS-0072154}
\date{June 15, 2001, revised September 12, 2001}
\maketitle
\begin{abstract}
In the first part of this paper,
given a smooth
family of Dirac-type operators on an odd-dimensional closed manifold,
we construct an abelian gerbe-with-connection whose curvature is the
three-form component of the Atiyah-Singer families index theorem.
In the second part of the paper, given a 
smooth family of Dirac-type operators
whose index lies in the subspace $\widetilde{\KK}^*_i(B)$ of the reduced
$\KK$-theory of the parametrizing space, 
we construct a set of Deligne cohomology classes of degree $i$
whose curvatures are the 
$i$-form component of the Atiyah-Singer families index theorem.
\end{abstract}

\section{Introduction}

To a family of $\overline{\partial}$-type operators  on a Hermitian
vector bundle over
a Riemann surface,  Quillen
associated the so-called determinant line bundle
\cite{Quillen (1985)}, a line bundle on the
parametrizing space with a natural connection. He also computed the
curvature of the connection. Quillen's construction was extended by Bismut and
Freed to the setting of a family of
Dirac-type operators on an even-dimensional closed 
manifold
\cite[Chapter 9.7]{Berline-Getzler-Vergne (1992)},
\cite{Bismut-Freed (1986)}.
The curvature of the connection is the two-form
component of the Atiyah-Singer families index theorem. A remarkable
feature of the determinant line bundle is that it is well-defined and
smooth even though the kernels and cokernels of the operators may not 
form vector bundles on the parametrizing space, due
to jumps in their dimensions.

In the first part of this paper 
we perform an analogous construction for a family of
Dirac-type operators on an odd-dimensional manifold. The determinant line
bundle is replaced by an abelian gerbe-with-connection. The ``curvature''
of the connection is the degree-$3$ component of the local families
index theorem, a $3$-form on the parametrizing space $B$.

In the second part of the paper we give a partial extension to the case of
degree $i \: > \: 3$.
Recall that the equivalence classes of line bundles with connection on 
$B$ are classified by the 2-dimensional Deligne cohomology of
$B$ \cite[Theorem 2.2.11]{Brylinski (1993)}. Similarly, the
equivalence classes of $\C^*$-gerbes-with-connection on $B$ 
are classified by the 3-dimensional Deligne cohomology of
$B$ \cite[Theorem 5.3.11]{Brylinski (1993)}. Hence for $i \: > \: 3$, 
in order to realize the
degree-$i$ component of the local families index theorem as the
``curvature'' of something, it is natural to
look for an $i$-dimensional Deligne cohomology class. 

There is an
apparent integrality obstruction to doing so, as
when $i \: > \: 3$ the degree-$i$ component of the
Chern character of the index class generally does not
lie in the image of the map $\HH^i(B; \Z) \rightarrow \HH^i(B; \Q)$. 
Hence we make an integrality assumption.
Recall that there is a filtration
${\KK}^*(B) = 
{\KK}^*_0(B) \supset {\KK}^*_1(B) \supset \ldots$ of 
the $\KK$-theory of $B$, where
${\KK}^*_i(B)$ consists of the elements $x$ of ${\KK}^*(B)$
with the property that for any finite simplicial complex $Y$ 
of dimension less than
$i$ and any continuous map $f \: : \: Y \rightarrow B$, 
$f^* x \: = \: 0$
\cite[Section 2]{Atiyah-Hirzebruch (1961)}. There is a similar filtration
of the reduced $\KK$-theory $\widetilde{\KK}^*(B)$. Taking $i$ to have the
parity of the dimension of the manifold on which the Dirac operator acts,
it turns out that we want to
assume that the image of the index of the family,
under the map ${\KK}^*(B) \rightarrow \widetilde{\KK}^*(B)$,
lies in $\widetilde{\KK}^*_i(B)$. Under this assumption,
we construct a countable set of explicit 
degree-$i$ Deligne cohomology classes on $B$ which only
depend on the geometrical input and whose ``curvatures'' are the
degree-$i$ component of the local families index theorem. 
(Note that if $i \: > \: 0$ and $\dim(B) \: > \: 0$ then 
the degree-$i$ Deligne cohomology of $B$ is infinite-dimensional.)
Roughly speaking, the different Deligne cohomology classes in the set 
correspond
to different trivializations of the index bundle on the $(i-2)$-skeleton
of a triangulation of $B$. Such trivializations are labeled by
$\bigoplus_{j=1}^\infty \HH^{i-1-2j}(B; \Z)$.

As a special case, if $B$ is $(i - 2)$-connected then 
the image of the index of the family automatically
lies in $\widetilde{\KK}^*_i(B)$, and we construct a unique
Deligne cohomology class.

In the rest of this introduction, we give an explicit statement of the
gerbe result. We defer the statement of the Deligne cohomology results
to Section \ref{Deligne} (see Theorem
\ref{thm2}). 

Information about gerbes is in the book of Brylinski
\cite{Brylinski (1993)} and the paper of Breen-Messing
\cite{Breen-Messing (2001)}. We will use a concrete approach to abelian
gerbes described by Hitchin \cite{Hitchin (1999)}. 
Given a manifold $B$ with a covering $\{U_\alpha\}_{\alpha \in I}$ by
open subsets, one obtains
a $\C^*$-gerbe from \\
1. A line bundle $L_{\alpha \beta}$ on each nonempty intersection
$U_\alpha \cap U_\beta$, \\
2. An isomorphism $L_{\alpha \beta} \: \cong \: L_{\beta \alpha}^{-1}$ and \\
3. A nowhere-zero section $\theta_{\alpha \beta \gamma}$ of
$L_{\alpha \beta} \: \otimes L_{\beta \gamma} \: \otimes L_{\gamma \alpha}$
on each nonempty intersection $U_\alpha \cap U_\beta \cap U_\gamma$ such
that \\
4. $\theta_{\beta \gamma \delta} \:
\theta_{\alpha \gamma \delta}^{-1} \:
\theta_{\alpha \beta \delta} \: 
\theta_{\alpha \beta \gamma}^{-1} \: = \: 1$ on 
each nonempty intersection $U_\alpha \cap U_\beta \cap U_\gamma \cap 
U_\delta$.

Given another choice 
$\left(L^\prime_{\alpha \beta}, \theta^\prime_{\alpha \beta
\gamma} \right)$, if there are line bundles 
$\{L_\alpha\}_{\alpha \in I}$ on the $U_\alpha$'s such
that $L^\prime_{\alpha \beta} \: \cong \: L_\alpha^{-1} \: \otimes \: 
L_{\alpha \beta} \:
\otimes \: L_{\beta}$, and $\theta^\prime_{\alpha \beta \gamma}$ is related
to $\theta_{\alpha \beta \gamma}$ in the obvious way,
then $\left(L^\prime_{\alpha \beta}, \theta^\prime_{\alpha \beta
\gamma} \right)$ is isomorphic to 
$\left(L_{\alpha \beta}, \theta_{\alpha \beta
\gamma} \right)$.  Taking a direct
limit over open coverings, one obtains the isomorphism classes of
gerbes on $B$. They are classified by
$\HH^3(B; \Z)$. 

A unitary connection on the gerbe
$\left(L_{\alpha \beta}, \theta_{\alpha \beta
\gamma} \right)$ is given by the additional data of \\
1. A unitary connection $\nabla_{\alpha \beta}$ on each $L_{\alpha \beta}$ 
and \\
2. A real $2$-form $F_\alpha \in \Omega^2(U_\alpha)$ on each $U_\alpha$\\
such that \\
1. $\nabla_{\alpha \beta} \: = \: \nabla_{\beta \alpha}^{-1}$,\\
2. $\theta_{\alpha \beta \gamma}$ is covariantly-constant with respect to
$\nabla_{\alpha \beta} \: \otimes \nabla_{\beta \gamma} \: \otimes 
\nabla_{\gamma \alpha}$ and \\
3. On each nonempty intersection $U_\alpha \cap U_\beta$, we have
$F_\beta \: - \: F_\alpha \: = \: 
c_1(\nabla_{\alpha \beta})$, the first Chern form of the connection
$\nabla_{\alpha \beta}$.

Suppose that 
$L^\prime_{\alpha \beta} \: = \: L_\alpha^{-1} \: \otimes \: 
L_{\alpha \beta} \:
\otimes \: L_{\beta}$ has connection
$\left( \nabla^\prime_{\alpha \beta}, F^\prime_\alpha \right)$.
If there are 
unitary connections $\nabla_\alpha$ on $L_\alpha$ such that \\
1. $\nabla^\prime_{\alpha \beta} \: = \: 
\nabla_\alpha^{-1} \: \otimes \: \nabla_{\alpha \beta}
\: \otimes \: \nabla_{\beta}$ and \\
2. $F^\prime_\alpha \: = \: F_\alpha 
+ c_1(\nabla_\alpha)$ \\
then
$\left( \nabla^\prime_{\alpha \beta}, F^\prime_\alpha \right)$ and
$\left( \nabla_{\alpha \beta}, F_\alpha \right)$
are equivalent .

The curvature of the connection, a globally-defined
$3$-form on $B$, is given on $U_\alpha$
by $d F_\alpha$.

Now let 
$\pi \: : \: M \rightarrow B$ be a smooth fiber bundle with closed
odd-dimensional fiber
$Z$. 
Let $TZ \: = \: \Ker(d\pi)$ denote the vertical tangent bundle, a
tangent bundle on $M$. We assume that $TZ$ has a spin structure. Let
$SZ$ be the corresponding spinor bundle.
Let $g^{TZ}$ be a vertical Riemannian metric. Let $V$ be a complex vector
bundle on $M$ with Hermitian metric $h^V$ and compatible
Hermitian connection $\nabla^V$.
Put $E \: = \: SZ \: \otimes \: V$. There is an ensuing
family $D_0 \: = \: \{(D_0)_b\}_{b \in B}$ of Dirac-type operators, with 
$(D_0)_b$ acting on
$C^\infty \left(Z_b; E \big|_{Z_b} \right)$. 

Let $T^HM$ be a horizontal
distribution on $M$. We now describe a gerbe on $B$.
We first choose an open covering $\{U_\alpha\}_{\alpha \in I}$ of $B$
with the property that there are
functions $\{h_\alpha\}_{\alpha \in I}$ in $C^\infty_c(\R)$ so that
$D_\alpha \: = \: D_0 \: + \: h_\alpha(D_0)$ is everywhere
invertible on $U_\alpha$.
It is easy to see that such $\{U_\alpha\}_{\alpha \in I}$ 
and $\{h_\alpha\}_{\alpha \in I}$ exist.  
If $U_\alpha \cap U_\beta \neq \emptyset$,
then the eigenvalues of the operators $
\frac{D_\beta}{|D_\beta|} \: - \: \frac{D_\alpha}{|D_\alpha|}$ 
over $U_\alpha \cap U_\beta$ are $0$, $2$ and $-2$. Let
$E_{-+}$ be orthogonal projection onto the eigenspace with eigenvalue $2$
and let 
$E_{+-}$ be orthogonal projection onto the eigenspace with eigenvalue $-2$.
(The notation for $E_{-+}$ is meant to indicate that on
$\Image(E_{-+})$, $D_\beta$ is positive and $D_\alpha$ is negative.) 
Then the images of $E_{-+}$ and $E_{+-}$ are finite-dimensional vector
bundles on $U_\alpha \cap U_\beta$.
Put
\begin{equation} \label{eqn1.1}
L_{\alpha \beta} \: = \: \Lambda^{max}(\Image(E_{-+})) \: \otimes
\left( \Lambda^{max}(\Image(E_{+-}) \right)^{-1}.
\end{equation}
If $U_\alpha \cap U_\beta \cap U_\gamma \: \neq \: \emptyset$ then
there is a canonical nowhere-zero section
$\Theta_{\alpha \beta \gamma}$ of 
$L_{\alpha \beta} \: \otimes L_{\beta \gamma} \: \otimes L_{\gamma \alpha}$
(see (\ref{eqn4.25})). 

The line bundle $L_{\alpha \beta}$ inherits a unitary connection
$\nabla_{\alpha \beta}$ from the projected connections on 
$\Image(E_{-+})$ and $\Image(E_{+-})$. We take
$F_\alpha$ to be the $2$-form component of a slight generalization of the
Bismut-Cheeger eta-form (see \cite[Definition 4.93]{Bismut-Cheeger (1989)} 
and (\ref{eqn2.19}) 
below). Usually in index theory the eta-form is most naturally
considered to be defined up to exact forms, but we will need the explicit
$2$-form component.

\begin{theorem} \label{thm1}
The data $\left( L_{\alpha \beta}, \theta_{\alpha \beta \gamma},
\nabla_{\alpha \beta}, F_\alpha \right)$ define a 
gerbe-with-connection on $B$ whose curvature is
\begin{equation} \label{eqn1.2}
\left( \int_Z \widehat{A} \left(R^{TZ}/2\pi i \right)
\: \wedge \: \ch \left( F^V/2\pi i \right) \right)^{(3)} \in \Omega^3(B).
\end{equation}
A different choice of $\{U_\alpha\}_{\alpha \in I}$ and 
$\{h_\alpha\}_{\alpha \in I}$ gives an equivalent
gerbe-with-connection.
\end{theorem}

Let us give a brief historical discussion of the relation between
gerbes and index theory.  This goes back to the index interpretation of
gauge anomalies.  Recall that from the Lagrangian viewpoint, the
nonabelian gauge anomaly arises from the possible
topological nontriviality of the
determinant line bundle on the space of connections modulo gauge
transformations \cite{Atiyah-Singer (1984)}.  
From the Hamiltonian viewpoint, this same anomaly 
becomes a $3$-dimensional cohomology class on the space of
connections modulo gauge transformations, namely the one that comes from the
families index theorem.  In \cite{Faddeev (1984)},
Faddeev constructs a $2$-cocycle on the gauge group which transgresses 
this $3$-dimensional cohomology class.  He interprets the cocycle as an
obstruction to satisfying Gauss' law.
In \cite[p. 200]{Pressley-Segal (1986)} Pressley and Segal note that
projective Hilbert bundles on $B$ are classified by $\HH^3(B; \Z)$, and
they use this to view the gauge anomaly as an obstruction to the
gauge-invariant construction of fermionic Fock spaces. Gerbes
(without connection) were brought into the picture by
Carey-Mickelsson-Murray \cite{Carey-Mickelsson-Murray (1997)},
Carey-Murray \cite{Carey-Murray (1996)} and Ekstrand-Mickelsson
\cite{Ekstrand-Mickelsson (2000)}.

Richard Melrose informs me that he and collaborators are working on 
related questions from a different viewpoint.  I thank Richard, 
Ulrich Bunke, Dan Freed,
Paolo Piazza, Stephan Stolz and Peter Teichner for discussions. 
I thank MSRI for its hospitality while this
research was performed.

\section{Conventions}

As for conventions, if $V$ is a vector bundle on $B$ with connection $\nabla^V$
and curvature $F^V \: = \: (\nabla^V)^2$ 
then we write $\ch(F^V) \: = \: \Tr \: \left( e^{- \: F^V} \right)
\in \Omega^{even}(B)$.
With this convention, $\ch(F^V/2\pi i) \: = \: \Tr \: \left( 
e^{- \: F^V/2 \pi i} \right)$ is
a closed form whose de Rham cohomology class lies
in the image of $\HH^*(B; \Q) \rightarrow \HH^*(B; \R)$.
We write $c_1(\nabla^V) \: = \: - \: \frac{1}{2 \pi i} \Tr(F^V) \: \in \:
\Omega^2(B)$. 

If $V$ is a $\Z_2$-graded
vector bundle on $B$ with a connection $\nabla^V$ that preserves the
$\Z_2$-grading,
and with curvature $F^V \: = \: (\nabla^V)^2$, 
then we write $\ch(F^V) \: = \: \Tr_s \: \left( e^{- \: F^V} \right)
\in \Omega^{even}(B)$.
Again,
$\ch(F^V/2\pi i) \: = \: \Tr_s \: \left( 
e^{- \: F^V/2 \pi i} \right)$ is
a closed form whose de Rham cohomology class lies
in the image of $\HH^*(B; \Q) \rightarrow \HH^*(B; \R)$.
We write $c_1(\nabla^V) \: = \: - \: \frac{1}{2 \pi i} \Tr_s(F^V) \: \in \:
\Omega^2(B)$.

If $g^{TB}$ is a
Riemannian metric on $B$ with curvature $2$-form $R^{TB}$ then we define
$\widehat{A}(R^{TB}) \in \Omega^{4*}(B)$ similarly, so that
$\widehat{A}(R^{TB}/2\pi i)$ is a de Rham representative of the usual
$\widehat{A}$-class in rational cohomology.

Let $\pi \: : \: M \rightarrow B$ be a smooth
fiber bundle as in the introduction, with fiber $Z$.
Let $T \in \Omega^2(M; TZ)$ denote the curvature of the horizontal
distribution, a $TZ$-valued horizontal $2$-form on $M$. 
Let $c(T)$ denote Clifford multiplication by $T$.

Let $\pi_* E$ be the infinite-dimensional vector
bundle on $B$ whose fiber over $b \in B$ is 
$C^\infty \left(Z_b; E \big|_{Z_b} \right)$. If $\dim(Z)$ is odd then
$\pi_* E$ is ungraded, while if $\dim(Z)$ is even then $\pi_* E$ is
$\Z_2$-graded.
Using $g^{TZ}$ and $h^V$, one obtains an $L^2$-inner
product $h^{\pi_* E}$ on $\pi_* E$. Let $\nabla^{\pi_* E}$ be
the canonical Hermitian connection on $\pi_* E$
\cite[Proposition 9.13]{Berline-Getzler-Vergne (1992)},
\cite[(4.21)]{Bismut-Cheeger (1989)}.

\section{The Index Gerbe}

\subsection{Eta-forms and their variations}

We now suppose that $Z$ is odd-dimensional.
Following \cite[\S 5]{Quillen (1985b)},
let $\sigma$ be a new formal odd variable such that 
$\sigma^2 \: = \: 1$. 

Let $D$ be the perturbation of $D_0$ by a smooth
family of fiberwise smoothing operators $P \: = \: \{P_b\}_{b \in B}$. That is,
$D_b \: = \: (D_0)_b \: + \: P_b$.
Given $s \: > \: 0$, the corresponding Bismut superconnection 
\cite[Section III]{Bismut (1986)}, 
\cite[Chapter 10.3]{Berline-Getzler-Vergne (1992)}
on 
$\pi_* E$ is
\begin{equation} \label{eqn2.1}
A_s \: = \: s \: \sigma \: D
\: + \: \nabla^{\pi_* E} \: + \: \frac{1}{4s} \: \sigma \: c(T).
\end{equation}
If $D \: = \: D_0$ then we write the superconnection as $A_{0,s}$.
Define $\Tr_\sigma$ on $(\C \: \oplus \: \C \sigma) \: \otimes \:
C^\infty(B; \End( \pi_* E))$ by
\begin{equation} \label{eqn2.2}
\Tr_\sigma (\alpha \: + \: \sigma \: \beta) \: = \:
\Tr(\beta) \in C^\infty(B),
\end{equation}
provided that $\alpha$ and $\beta$ are fiberwise
trace-class operators.  Then there
is an extension of $\Tr_\sigma$ to an
$\Omega^*(B)$-valued trace on
$\left( \Omega^*(B) \: \widehat{\otimes} \: 
(\C \: \oplus \: \C \sigma) \right) \: 
\otimes_{C^\infty(B)} C^\infty(B; \End( \pi_* E))$ which is
left-$\Omega^*(B)$ linear,
again provided that the vertical
operators are trace-class.
For any $s \: > \: 0$,
$\Tr_\sigma \left( e^{- \: A_s^2} \right) \in \Omega^{odd}(B)$
represents the Chern character of the index $\Ind(D) \in \KK^1(B)$ of
the family of vertical operators, up
to normalizing constants.
For later use, we note that
\begin{align} \label{eqn2.3}
A_s^2 \: = \: s^2 \: D^2 \: & - \: s \: \sigma \: [\nabla, D] \: + \:
\left( \nabla^2 \: + \: \frac{1}{4} \: (D \: c(T) \: + \: c(T) \: D) \right) \\
& - \: \frac{1}{4s} \: \sigma \: [\nabla, c(T)] \: + \: 
\frac{1}{16 s^2} \: c(T)^2. \notag
\end{align}
The meaning of $\Tr_\sigma \left( e^{- \: A_s^2} \right)$ is that
the component in $\Omega^{2k+1}(B)$ is derived by means of a Duhamel expansion
around $e^{- \: s^2 \: D^2}$
\cite[Appendix to Chapter 9]{Berline-Getzler-Vergne (1992)}, 
and hence comes from a finite number of terms
in the Duhamel expansion.

If $D \: = \: D_0$ then
$\lim_{s \rightarrow 0} \Tr_\sigma \left( e^{- \: A_{0,s}^2} \right)$ exists
and \cite[(4.97)]{Bismut-Cheeger (1989)}
\begin{equation} \label{eqn2.4}
\lim_{s \rightarrow 0} \Tr_\sigma \left( e^{- \: A_{0,s}^2} \right) \: = \:
\sqrt{\pi} \: (2 \pi i)^{- \:
\frac{dim(Z) \: + \: 1}{2}} \:
\int_Z \widehat{A} \left(R^{TZ} \right)
\: \wedge \: \ch \left( F^V \right).
\end{equation}
(The constants in this expression will most easily be seen as arising
from (\ref{eqn2.26})).
For general $D$, we do not know that
$\Tr_\sigma \left( e^{- \: A_{s}^2} \right)$ has a limit as
$s \rightarrow 0$. However, let $\LIM_{s \rightarrow 0}$ denote the principal 
value as in \cite[Section 9.6]{Berline-Getzler-Vergne (1992)}. Then
by expanding in a Duhamel series around 
$\Tr_\sigma \left( e^{- \: A_{0,s}^2} \right)$, one finds that
$\LIM_{s \rightarrow 0} \Tr_\sigma \left( e^{- \: A_s^2} \right)$ 
exists.

\begin{proposition} \label{prop1}
For all $D$,
\begin{equation} \label{eqn2.5}
\LIM_{s \rightarrow 0} \:
\Tr_\sigma \left( e^{- \: A_s^2} \right) 
\: = \: \sqrt{\pi} \: (2 \pi i)^{- \:
\frac{dim(Z) \: + \: 1}{2}} \:
\int_Z \widehat{A} \left( R^{TZ} \right)
\: \wedge \: \ch \left( F^V \right).
\end{equation}
\end{proposition}
\begin{pf}
In general, if $\{A_s(\epsilon)\}_{\epsilon \in [0,1]}$ is a smooth
$1$-parameter family of superconnections then formally,
\begin{align} \label{eqn2.6}
\frac{d}{d\epsilon} \Tr_\sigma \left( e^{- \: A_s(\epsilon)^2} \right) \: & =
- \: \Tr_\sigma \left( \left\{ A_s(\epsilon), 
\frac{dA_s(\epsilon)}{d\epsilon} \right\} \:
e^{- \: A_s(\epsilon)^2} \right) \\
& =
- \: \Tr_\sigma \left( \left\{ A_s(\epsilon),
\frac{dA_s(\epsilon)}{d\epsilon} \:
e^{- \: A_s(\epsilon)^2} \right\} \right) \notag \\
& = \: - \: d \: \Tr_\sigma \left( \frac{dA_s(\epsilon)}{d\epsilon} \:
e^{- \: A_s(\epsilon)^2} \right). \notag
\end{align}
Let $\{D(\epsilon)\}_{\epsilon \in [0,1]}$ be a smooth $1$-parameter family
of operators $D$ as above. Then (\ref{eqn2.6}) is easily justified, and gives
\begin{equation} \label{eqn2.7}
\frac{d}{d\epsilon} \Tr_\sigma \left( e^{- \: A_s(\epsilon)^2} \right) \: = \:
- \: d \: \Tr_\sigma \left( s \: \sigma \: \frac{dD(\epsilon)}{d\epsilon} \:
e^{- \: A_s(\epsilon)^2} \right).
\end{equation}
Hence
\begin{equation} \label{eqn2.8}
\LIM_{s \rightarrow 0} \:
\frac{d}{d\epsilon} \Tr_\sigma \left( e^{- \: A_s(\epsilon)^2} \right) \: = \:
- \: d \: 
\LIM_{s \rightarrow 0} \:
 \Tr_\sigma \left( s \: \sigma \: \frac{dD(\epsilon)}{d\epsilon} \:
e^{- \: A_s(\epsilon)^2} \right).
\end{equation}
The Duhamel expansion of $
s \: \sigma \: \frac{dD(\epsilon)}{d\epsilon} \:
e^{- \: A_s(\epsilon)^2}$ is
\begin{align} \label{eqn2.9}
\sum_{l=0}^\infty (-1)^l \: \int_0^1 \ldots \int_0^1 \:
& s \: \sigma \: \frac{dD}{d\epsilon} \: e^{- \: t_0 \:
s^2 \: D^2} \: (A_s(\epsilon)^2 \: - \: s^2 \: D^2) \: e^{- \: t_1 \:
s^2 \: D^2} \: (A_s(\epsilon)^2 \: - \: s^2 \: D^2) \ldots \\
& (A_s(\epsilon)^2 
\: - \: s^2 \: D^2) \:
e^{- \: t_l \: s^2 \: D^2} 
\: \delta(t_0 \: + \: \ldots \: + t_l \: - \: 1)
\: dt_0 \ldots dt_l. \notag
\end{align}
If we consider the component of (\ref{eqn2.9})
of degree $2k$ with respect to $B$
then only a finite number of terms in the expansion
(\ref{eqn2.9}) will enter. From (\ref{eqn2.3}), 
\begin{equation} \label{eqn2.10}
A_s^2 \: - \: s^2 \: D^2 \: = \: s^{-2} f(s \sigma)
\end{equation}
for a polynomial $f$ with appropriate coefficients.  
As $\frac{dD}{d\epsilon}$ is smoothing, we can compute 
$\LIM_{s \rightarrow 0} \:
\Tr_\sigma \left( s \: \sigma \: \frac{dD(\epsilon)}{d\epsilon} \:
e^{- \: A_s(\epsilon)^2} \right)$ by looking at the terms of (\ref{eqn2.9}) 
which contribute to $\Omega^{2k}(B)$ and expanding the exponentials
$e^{- \: t_j \: s^2 \: D^2}$ in $s^2$. In so doing, 
the resulting expression is a Laurent series of the form
$s \: \sigma \: s^{- \: 2L} \: \sum_{r=0}^\infty c_r \: (s \sigma)^r$ 
for some $L \: \ge \: 0$.
Then after applying $\Tr_\sigma$, the result
is $s^{1 \: - \: 2L} \: \sum_{r \: even} c_r \: s^r$.
Hence 
\begin{equation} \label{eqn2.11}
\LIM_{s \rightarrow 0} \: 
\Tr_\sigma \left( 
s \: \sigma \: \frac{d D(\epsilon)}{d \epsilon} \: 
e^{- \: A_s(\epsilon)^2} \right) \: = \: 0
\end{equation}
and so from (\ref{eqn2.8}),
\begin{equation} \label{eqn2.12}
\LIM_{s \rightarrow 0} \:
\frac{d}{d\epsilon} \Tr_\sigma \left( e^{- \: A_s(\epsilon)^2} \right) \: = \:
0.
\end{equation} 

In our case, we can commute $\LIM_{s \rightarrow 0}$ and
$\frac{d}{d\epsilon}$. Taking 
\begin{equation} \label{eqn2.13}
D(\epsilon) \: = \: D_0 \: + \: \epsilon \: (D \: - \: D_0),
\end{equation}
we obtain
\begin{align} \label{eqn2.14}
\LIM_{s \rightarrow 0} \:
\Tr_\sigma \left( e^{- \: A_s^2} \right) \: & = \:
\LIM_{s \rightarrow 0} \:
\Tr_\sigma \left( e^{- \: A_{0,s}^2} \right) \: = \:
\lim_{s \rightarrow 0} \:
\Tr_\sigma \left( e^{- \: A_{0,s}^2} \right) \\
& = \: \sqrt{\pi} \: (2 \pi i)^{- \:
\frac{dim(Z) \: + \: 1}{2}} \:
\int_Z \widehat{A} \left( R^{TZ} \right)
\: \wedge \: \ch \left( F^V \right). \notag
\end{align} 
\end{pf}

The next result is the same as \cite[Proposition 14]{Melrose-Piazza (1997)}.
We give the proof for completeness.
\begin{proposition} \label{prop2}
If $\{D(\epsilon)\}_{\epsilon \in [0,1]}$ is a smooth $1$-parameter family
of operators as above then
\begin{align} \label{eqn2.15}
\frac{\partial}{\partial \epsilon} \Tr_\sigma 
\left( \frac{\partial A_s(\epsilon)}{\partial s} \: 
e^{- \: A_s(\epsilon)^2} \right) \: & - \:
\frac{\partial}{\partial s} \Tr_\sigma \left( 
\frac{\partial A_s(\epsilon)}{\partial \epsilon} \: 
e^{- \: A_s(\epsilon)^2} \right) \: = \\
& d \int_0^1 \Tr_\sigma \left(  
\frac{\partial A_s(\epsilon)}{\partial \epsilon}
\: e^{- \: u \: A_s(\epsilon)^2} \: 
\frac{\partial A_s(\epsilon)}{\partial s} \:
e^{- \: (1-u) \: A_s(\epsilon)^2} \right) \: du. \notag
\end{align}
\end{proposition}
\begin{pf}
Formally,
\begin{align} \label{eqn2.16}
\frac{\partial}{\partial \epsilon} \Tr_\sigma 
\left( \frac{\partial A_s(\epsilon)}{\partial s} \: 
e^{- \: A_s(\epsilon)^2} \right) \: = \: 
& \Tr_\sigma 
\left( \frac{\partial^2 A_s(\epsilon)}{\partial \epsilon \: \partial s} \: 
e^{- \: A_s(\epsilon)^2} \right) \: - \\
& \int_0^1 \Tr_\sigma \left(  
\frac{\partial A_s(\epsilon)}{\partial \epsilon}
\: e^{- \: u \: A_s(\epsilon)^2} \: 
\left\{ A_s(\epsilon), \frac{\partial A_s(\epsilon)}{\partial s} \right\} \:
e^{- \: (1-u) \: A_s(\epsilon)^2} \right) \: du \notag
\end{align}
and
\begin{align} \label{eqn2.17}
\frac{\partial}{\partial s} \Tr_\sigma 
\left( \frac{\partial A_s(\epsilon)}{\partial s} \: 
e^{- \: A_s(\epsilon)^2} \right) \: = \: 
& \Tr_\sigma 
\left( \frac{\partial^2 A_s(\epsilon)}{\partial s \: \partial \epsilon} \: 
e^{- \: A_s(\epsilon)^2} \right) \: - \\
& \int_0^1 \Tr_\sigma \left(  
\frac{\partial A_s(\epsilon)}{\partial s}
\: e^{- \: u \: A_s(\epsilon)^2} \: 
\left\{ A_s(\epsilon), \frac{\partial A_s(\epsilon)}{\partial 
\epsilon} \right\} \:
e^{- \: (1-u) \: A_s(\epsilon)^2} \right) \: du. \notag
\end{align}
Then
\begin{align} \label{eqn2.18}
& \frac{\partial}{\partial \epsilon} \Tr_\sigma 
\left( \frac{\partial A_s(\epsilon)}{\partial s} \: 
e^{- \: A_s(\epsilon)^2} \right) \: - \: 
\frac{\partial}{\partial s} \Tr_\sigma 
\left( \frac{\partial A_s(\epsilon)}{\partial s} \: 
e^{- \: A_s(\epsilon)^2} \right) \: = \\
& \int_0^1 \Tr_\sigma \left(  
\frac{\partial A_s(\epsilon)}{\partial s}
\: e^{- \: u \: A_s(\epsilon)^2} \: 
\left\{ A_s(\epsilon), \frac{\partial A_s(\epsilon)}{\partial 
\epsilon} \right\} \:
e^{- \: (1-u) \: A_s(\epsilon)^2} \right) \: du \: - \notag \\
& \int_0^1 \Tr_\sigma \left(  
\frac{\partial A_s(\epsilon)}{\partial \epsilon}
\: e^{- \: u \: A_s(\epsilon)^2} \: 
\left\{ A_s(\epsilon), \frac{\partial A_s(\epsilon)}{\partial s} \right\} \:
e^{- \: (1-u) \: A_s(\epsilon)^2} \right) \: du \: = \: \notag \\
& \int_0^1 \Tr_\sigma \left(  
\left\{ A_s(\epsilon), \frac{\partial A_s(\epsilon)}{\partial 
\epsilon} \right\} \:
e^{- \: (1-u) \: A_s(\epsilon)^2} \: 
\frac{\partial A_s(\epsilon)}{\partial s}
\: e^{- \: u \: A_s(\epsilon)^2} \: 
 \right) \: du \: - \notag \\
& \int_0^1 \Tr_\sigma \left(  
\frac{\partial A_s(\epsilon)}{\partial \epsilon}
\: e^{- \: u \: A_s(\epsilon)^2} \: 
\left\{ A_s(\epsilon), \frac{\partial A_s(\epsilon)}{\partial s} \right\} \:
e^{- \: (1-u) \: A_s(\epsilon)^2} \right) \: du \: = \: \notag \\
& \int_0^1 \Tr_\sigma \left(  
\left\{ A_s(\epsilon), \frac{\partial A_s(\epsilon)}{\partial 
\epsilon} \right\} \:
e^{- \: u \: A_s(\epsilon)^2} \: 
\frac{\partial A_s(\epsilon)}{\partial s}
\: e^{- \: (1 - u) \: A_s(\epsilon)^2} \: 
 \right) \: du \: - \notag \\
& \int_0^1 \Tr_\sigma \left(  
\frac{\partial A_s(\epsilon)}{\partial \epsilon}
\: e^{- \: u \: A_s(\epsilon)^2} \: 
\left\{ A_s(\epsilon), \frac{\partial A_s(\epsilon)}{\partial s} \right\} \:
e^{- \: (1-u) \: A_s(\epsilon)^2} \right) \: du \: = \: \notag \\
& \int_0^1 \Tr_\sigma \left(  
\left\{ A_s(\epsilon), \frac{\partial A_s(\epsilon)}{\partial 
\epsilon} \:
e^{- \: u \: A_s(\epsilon)^2} \: 
\frac{\partial A_s(\epsilon)}{\partial s}
\: e^{- \: (1 - u) \: A_s(\epsilon)^2} \: du  \right\} \:
 \right) \: = \notag \\
&\int_0^1 d \: \Tr_\sigma \left(  
\frac{\partial A_s(\epsilon)}{\partial 
\epsilon} \:
e^{- \: u \: A_s(\epsilon)^2} \: 
\frac{\partial A_s(\epsilon)}{\partial s}
\: e^{- \: (1 - u) \: A_s(\epsilon)^2} \: 
 \right) \: du \: = \: \notag \\
& d \: \int_0^1 \Tr_\sigma \left(  
\frac{\partial A_s(\epsilon)}{\partial 
\epsilon} \:
e^{- \: u \: A_s(\epsilon)^2} \: 
\frac{\partial A_s(\epsilon)}{\partial s}
\: e^{- \: (1 - u) \: A_s(\epsilon)^2} \: 
 \right) \: du. \notag
\end{align}
It is easy to justify these formal manipulations.
\end{pf}

Let $\{U_\alpha\}_{\alpha \in I}$ 
be a covering of $B$ by open sets. For each $\alpha \in I$,
suppose that $D_\alpha$ is a family of operators 
$\{(D_\alpha)_b\}_{b \in U_\alpha}$ as before, defined over $U_\alpha$. 
We assume that
for each $b \in U_\alpha$, $(D_\alpha)_b$ is invertible.
Given $\{(D_0)_b\}_{b \in U_\alpha}$, 
the obstruction to finding such a family $D_\alpha$ is the index of
$D_0 \big|_{U_\alpha}$ in 
$\KK^1(U_\alpha)$ \cite[Proposition 1]{Melrose-Piazza (1997)}. 
For example, if $U_\alpha$ is contractible then 
there is no obstruction.

From the method of proof of 
\cite[Theorem 10.32]{Berline-Getzler-Vergne (1992)},
$\Tr_\sigma \left( 
\frac{d A_{0,s}}{d s} \: 
e^{- \: A_{0,s}^2} \right)$ has an asymptotic expansion as $s \rightarrow 0$ of
the form $\sum_{k=0}^\infty a_k \: s^k$. Then by the method of proof of
Proposition \ref{prop1}, the degree-$2k$ component of
$\Tr_\sigma \left( 
\frac{d A_s}{d s} \: 
e^{- \: A_s^2} \right)$ will have an asymptotic expansion as
$s \rightarrow 0$ of the
form $s^{-2L} \: \sum_{k=0}^\infty b_k \: s^k$.  Hence from
\cite[Lemma 9.34]{Berline-Getzler-Vergne (1992)}, it makes sense to define
$\widetilde{\eta}_\alpha
\in \Omega^{even}(U_\alpha)$ by 
\begin{equation} \label{eqn2.19}
\widetilde{\eta}_\alpha
\: = \: \LIM_{t \rightarrow 0} \int_t^\infty
\Tr_\sigma \left( 
\frac{d A_s}{d s} \: 
e^{- \: A_s^2} \right) \: ds;
\end{equation}
compare \cite[Definition 2.4]{Dai-Zhang (1998)}.
As
\begin{equation} \label{eqn2.20}
\frac{d}{ds} \Tr_\sigma \left( e^{- \: A_s^2} \right) \: = \:
- \: d \: \Tr_\sigma \left( 
\frac{\partial A_s}{\partial s} \: 
e^{- \: A_s^2} \right),
\end{equation}
it follows that
\begin{equation} \label{eqn2.21}
d \widetilde{\eta}_\alpha \: = \:
\LIM_{s \rightarrow 0} \: \Tr_\sigma \left( e^{- \: A_s^2} \right) \: = \:
\sqrt{\pi} \: (2 \pi i)^{- \:
\frac{dim(Z) \: + \: 1}{2}} \:
\int_Z \widehat{A} \left( R^{TZ} \right)
\: \wedge \: \ch \left( F^V \right).
\end{equation}

Let $A_s(\epsilon)$ be a smooth $1$-parameter family of superconnections.
As in Proposition \ref{prop2}, when the terms make sense, we have
\begin{align} \label{eqn2.22}
\frac{d \widetilde{\eta}_\alpha(\epsilon)}{d\epsilon}  \: = \:
& - \: \LIM_{s \rightarrow 0} \: 
\Tr_\sigma \left( 
\frac{\partial A_s(\epsilon)}{\partial \epsilon} \: 
e^{- \: A_s(\epsilon)^2} \right) \: + \: \\
& d \: \LIM_{t \rightarrow 0} \int_t^\infty \int_0^1
\Tr_\sigma \left(  
\frac{\partial A_s(\epsilon)}{\partial \epsilon}
\: e^{- \: u \: A_s(\epsilon)^2} \: 
\frac{\partial A_s(\epsilon)}{\partial s} \:
e^{- \: (1-u) \: A_s(\epsilon)^2} \right) \: du \: ds. \notag
\end{align}

\begin{proposition} \label{prop3}
Let $\{D_\alpha(\epsilon)\}_{\epsilon \in [0,1]}$ be a smooth
$1$-parameter family of $D_\alpha$'s as before. Then
\begin{equation} \label{eqn2.23}
\frac{d \widetilde{\eta}_\alpha(\epsilon)}{d\epsilon}  \: = \:
d \: \LIM_{t \rightarrow 0} \int_t^\infty \int_0^1
\Tr_\sigma \left(  
\frac{\partial A_s(\epsilon)}{\partial \epsilon}
\: e^{- \: u \: A_s(\epsilon)^2} \: 
\frac{\partial A_s(\epsilon)}{\partial s} \:
e^{- \: (1-u) \: A_s(\epsilon)^2} \right) \: du \: ds.
\end{equation}
\end{proposition}
\begin{pf}
This follows from (\ref{eqn2.11}), Proposition \ref{prop2}
and  (\ref{eqn2.22}).
\end{pf}

Proposition \ref{prop3} is closely related to
\cite[Corollary 4]{Melrose-Piazza (1997)}.

To give normalizations that are compatible with rational cohomology,
let ${\cal R}$ be the operator on $\Omega^*(B)$ which acts on
$\Omega^{2k}(B)$ as multiplication by $(2\pi i)^{-k}$ and which acts
on $\Omega^{2k+1}(B)$ as multiplication by $(2\pi i)^{-k}$.
Put 
\begin{equation} \label{eqn2.24}
\ch(A_s) \: = \: \pi^{- \: 1/2} \: {\cal R} \: 
\Tr_\sigma \left( e^{- \: A_s^2} \right)
\end{equation}
and 
\begin{equation} \label{eqn2.25}
\widehat{\eta} \: = \: \pi^{- \: 1/2} \: {\cal R} \: \widetilde{\eta}.
\end{equation}
Then (\ref{eqn2.21}) becomes
\begin{equation} \label{eqn2.26}
d \widehat{\eta}_\alpha \: = \:
\LIM_{s \rightarrow 0} \: \ch(A_s) \: = \:
\int_Z \widehat{A} \left(R^{TZ}/2\pi i \right)
\: \wedge \: \ch \left( F^V/2\pi i \right).
\end{equation}

\subsection{The $1$-form case}

The degree-$0$ component $\widehat{\eta}_\alpha^{(0)} \in
\Omega^0(U_\alpha)$ of $\widehat{\eta}_\alpha$,
i.e.
\begin{align} \label{eqn3.1}
\widehat{\eta}_\alpha^{(0)} \: & = \:
\pi^{- \: 1/2} \: \LIM_{t \rightarrow 0} \: \int_t^\infty
\Tr_\sigma \left( \sigma \: D_\alpha \: e^{- \: s^2 \: D_\alpha^2} \right)
\: ds \: \\
& = \: 
\pi^{- \: 1/2} \: \LIM_{t \rightarrow 0} \: \int_t^\infty
\Tr \left(D_\alpha \: e^{- \: s^2 \: D_\alpha^2} \right) \: ds, \notag
\end{align}
is half of the
Atiyah-Singer-Patodi eta-invariant of $D_\alpha$
\cite{Atiyah-Patodi-Singer (1975)}. 
If $U_\alpha \cap U_\beta \: \neq \: \emptyset$ then it is well-known that
$\widehat{\eta}_\beta^{(0)} \big|_{U_\alpha \cap U_\beta} \: - \: 
\widehat{\eta}_\alpha^{(0)} \big|_{U_\alpha \cap U_\beta}$ is an
integer-valued function on $U_\alpha \cap U_\beta$. 
Hence if $f_\alpha \: :\: U_\alpha \rightarrow S^1$ is defined by
$f_\alpha \: = \: e^{2 \pi i \widehat{\eta}_\alpha^{(0)}}$ then
if $U_\alpha \cap U_\beta \: \neq \: \emptyset$,
$f_\alpha \big|_{U_\alpha \cap U_\beta} \: = \: 
f_\beta \big|_{U_\alpha \cap U_\beta}$. Thus the functions
$\{f_\alpha\}_{\alpha \in I}$ piece together to give a function
$f \: : \: B \rightarrow S^1$ such that
$f \big|_{U_\alpha} \: = \: f_\alpha$. 
From (\ref{eqn2.26}),
\begin{equation} \label{eqn3.2}
\frac{1}{2 \pi i} \: d \ln f \: = \: 
\left( \int_Z \widehat{A} \left(R^{TZ}/2\pi i \right)
\: \wedge \: \ch \left( F^V/2\pi i \right) \right)^{(1)} \in
\Omega^1(B).
\end{equation}
In particular, if
$[S^1] \in \HH^1(S^1; \Z)$ is the fundamental class of $S^1$ then
$f^* [S^1] \in \HH^1(B; \Z)$ is represented in real cohomology by the 
closed form on the right-hand-side of (\ref{eqn3.2}).

\subsection{The index gerbe}

Let $\widetilde{\eta}_\alpha^{(2)} \in
\Omega^2(U_\alpha)$ denote
the degree-$2$ component of $\widetilde{\eta}_\alpha$.

\begin{proposition} \label{prop4}
If $\{D_\alpha(\epsilon)\}_{\epsilon \in [0,1]}$ is a smooth
$1$-parameter family of mutually-commuting invertible operators as before then
$\widetilde{\eta}_\alpha(\epsilon)^{(2)} \in
\Omega^2(U_\alpha)$ is independent of $\epsilon$.
\end{proposition}
\begin{pf}
From Proposition \ref{prop3}, it is enough to show the vanishing of
the component of 
\begin{equation} \label{eqn4.1}
\LIM_{t \rightarrow 0} \int_t^\infty \int_0^1
\Tr_\sigma \left(  
\frac{\partial A_s(\epsilon)}{\partial \epsilon}
\: e^{- \: u \: A_s(\epsilon)^2} \: 
\frac{\partial A_s(\epsilon)}{\partial s} \:
e^{- \: (1-u) \: A_s(\epsilon)^2} \right) \: du \: ds
\end{equation}
in $\Omega^1(B)$.
This is the degree-$1$ component of
\begin{equation} \label{eqn4.2}
\LIM_{t \rightarrow 0} \int_t^\infty \int_0^1
\Tr_\sigma \left(  
s \: \sigma \: \frac{dD_\alpha}{d\epsilon}
\: e^{- \: u \: (s^2 \: D_\alpha^2 \: - \: s \: \sigma \: [\nabla, D_\alpha])}
\: 
\sigma \: D_\alpha \:
e^{- \: (1-u) \: 
(s^2 \: D_\alpha^2 \: - \: s \: \sigma \: [\nabla, D_\alpha])} 
\right) \: du \: ds.
\end{equation}
Using the fact that $\frac{dD_\alpha}{d\epsilon}$ commutes with
$D_\alpha$, the degree-$1$ component of 
\begin{equation} \label{eqn4.3}
\Tr_\sigma \left(  
s \: \sigma \: \frac{dD_\alpha}{d\epsilon}
\: e^{- \: u \: (s^2 \: D_\alpha^2 \: - \: s \: \sigma \: [\nabla, D_\alpha])}
\: 
\sigma \: D_\alpha \:
e^{- \: (1-u) \: 
(s^2 \: D_\alpha^2 \: - \: s \: \sigma \: [\nabla, D_\alpha])} 
\right)
\end{equation}
is
\begin{align} \label{eqn4.4}
& \Tr_\sigma \left( 
s \: \sigma \: \frac{dD_\alpha}{d\epsilon}
\:
u \: s \: \sigma \: [\nabla, D_\alpha]
\: e^{- \: u \: s^2 \: D_\alpha^2}
\: 
\sigma \: D_\alpha \:
e^{- \: (1-u) \: 
s^2 \: D_\alpha^2}
\: + \: \right. \\
& \left. s \: \sigma \: \frac{dD_\alpha}{d\epsilon}
\: e^{- \: u \: s^2 \: D_\alpha^2}
\: 
\sigma \: D_\alpha \:
e^{- \: (1-u) \: 
s^2 \: D_\alpha^2}
\: (1-u) \: 
s \: \sigma \: [\nabla, D_\alpha]
\right) \: = \notag \\
& \Tr_\sigma \left(  
s \: \sigma \: \frac{dD_\alpha}{d\epsilon}
\:
u \: s \: \sigma \: [\nabla, D_\alpha]
\: e^{- \: s^2 \: D_\alpha^2}
\: 
\sigma \: D_\alpha
\: + \:
\right. \notag
\\ & \left. 
s \: \sigma \: \frac{dD_\alpha}{d\epsilon}
\: e^{- \: s^2 \: D_\alpha^2}
\: 
\sigma \: D_\alpha \:
\: (1-u) \: 
s \: \sigma \: [\nabla, D_\alpha]
\right) \: = \notag \\
& \Tr_\sigma \left(  - \: 
s \: \sigma \: \frac{dD_\alpha}{d\epsilon}
\:
u \: s \: \sigma
\: e^{- \: s^2 \: D_\alpha^2}
\: 
\sigma \: D_\alpha \:  [\nabla, D_\alpha]
\: + \:
\right. \notag \\ & \left. 
s \: \sigma \: \frac{dD_\alpha}{d\epsilon}
\: e^{- \: s^2 \: D_\alpha^2}
\: 
\sigma \: D_\alpha \:
\: (1-u) \: 
s \: \sigma \: [\nabla, D_\alpha]
\right) \: = \notag \\
& (1 \: - \: 2 \: u) \: s^2 \: \Tr \left( 
\frac{dD_\alpha}{d\epsilon}
\: e^{- \: s^2 \: D_\alpha^2}
\: D_\alpha \:  [\nabla, D_\alpha] \right). \notag
\end{align}
As $\int_0^1 (1 \: - \: 2 \: u) \: du \: = \: 0$, the proposition follows.
\end{pf}

\subsection{Finite-dimensional case}

Let $V$ be a finite-dimensional Hermitian vector bundle with a connection
$\nabla$ over $B$.  Let $D \in \End(V)$ be an invertible self-adjoint
operator. 
As in \cite[Section 2(b)]{Bismut-Cheeger (1989)},
put $A_s \: = \: s \: \sigma \: D \: + \: \nabla$ and
\begin{equation} \label{eqn4.5}
\widetilde{\eta}
\: = \: \int_0^\infty
\Tr_\sigma \left( 
\frac{d A_s}{d s} \: 
e^{- \: A_s^2} \right) \: ds \in \Omega^{even}(B).
\end{equation}
Put $P_\pm \: = \: \frac{|D| \: \pm \: D}{2 |D|}$.
From 
\cite[Theorem 2.43]{Bismut-Cheeger (1989)}, $\widetilde{\eta}$ is closed and,
up to normalizing constants,
represents the Chern character of $[\Image(P_+) \: - \: \Image(P_-)] \in
\KK_0(B)$ in $\HH^{even}(B; \C)$. We wish to say precisely what
$\widetilde{\eta}^{(2)} \in \Omega^2(B)$ is.

\begin{proposition} \label{prop5}
\begin{equation} \label{eqn4.6}
\widetilde{\eta}^{(2)} \: = \: - \: \frac{\sqrt{\pi}}{2} \:
\Tr \left( \left( P_+ \: \nabla \:
P_+ \right)^2 \: - \:  
\left( P_- \: \nabla \:
P_- \right)^2 \right) \in \Omega^2(B).
\end{equation}
\end{proposition}
\begin{pf}
Let $\{\nabla(\epsilon)\}_{\epsilon \in [0,1]}$ be a smooth
$1$-parameter family of Hermitian connections on $V$. As in 
(\ref{eqn2.22}), we have
\begin{align} \label{eqn4.7}
\frac{d \widetilde{\eta}(\epsilon)}{d\epsilon}  \: & = \:
 - \: \lim_{s \rightarrow 0} \: 
\Tr_\sigma \left( 
\frac{\partial A_s(\epsilon)}{\partial \epsilon} \: 
e^{- \: A_s(\epsilon)^2} \right) \: + \\
& \: \: \: \: \: \: \: d \: \int_0^\infty \int_0^1
\Tr_\sigma \left(  
\frac{\partial A_s(\epsilon)}{\partial \epsilon}
\: e^{- \: u \: A_s(\epsilon)^2} \: 
\frac{\partial A_s(\epsilon)}{\partial s} \:
e^{- \: (1-u) \: A_s(\epsilon)^2} \right) \: du \: ds  \notag \\
& = \:
 - \: \lim_{s \rightarrow 0} \: 
\Tr_\sigma \left( 
\frac{d \nabla(\epsilon)}{d \epsilon} \: 
e^{- \: A_s(\epsilon)^2} \right) \: + \notag \\
& \: \: \: \: \: \: \: d \: \int_0^\infty \int_0^1
\Tr_\sigma \left(  
\frac{d\nabla(\epsilon)}{d \epsilon}
\: e^{- \: u \: A_s(\epsilon)^2} \: 
\sigma \: D \:
e^{- \: (1-u) \: A_s(\epsilon)^2} \right) \: du \: ds \notag \\
& = \: - \:
\Tr_\sigma \left( 
\frac{d \nabla(\epsilon)}{d \epsilon} \: 
e^{- \: \nabla(\epsilon)^2} \right) \: + \notag \\
& \: \: \: \: \: \: \: d \: \int_0^\infty \int_0^1
\Tr_\sigma \left(  
\frac{d\nabla(\epsilon)}{d \epsilon}
\: e^{- \: u \: A_s(\epsilon)^2} \: 
\sigma \: D \:
e^{- \: (1-u) \: A_s(\epsilon)^2} \right) \: du \: ds \notag \\
& = \: d \: \int_0^\infty \int_0^1
\Tr_\sigma \left(  
\frac{d\nabla(\epsilon)}{d \epsilon}
\: e^{- \: u \: A_s(\epsilon)^2} \: 
\sigma \: D \:
e^{- \: (1-u) \: A_s(\epsilon)^2} \right) \: du \: ds. \notag
\end{align}
The $1$-form component of 
$\Tr_\sigma \left(  
\frac{d\nabla(\epsilon)}{d \epsilon}
\: e^{- \: u \: A_s(\epsilon)^2} \: 
\sigma \: D \:
e^{- \: (1-u) \: A_s(\epsilon)^2} \right)$ is
\begin{equation} \label{eqn4.8}
\Tr_\sigma \left(  
\frac{d\nabla(\epsilon)}{d \epsilon}
\: e^{- \: u \: s^2 \: D^2} \: 
\sigma \: D \:
e^{- \: (1-u) \: s^2 \: D^2} \right) \: = \: - \:
\Tr \left(  
\frac{d\nabla(\epsilon)}{d \epsilon}
\: D \:
e^{- \: s^2 \: D^2} \right).
\end{equation}
Then
\begin{equation} \label{eqn4.9}
\frac{d \widetilde{\eta}(\epsilon)^{(2)}}{d\epsilon}  \: = \: - \:
d \: \int_0^\infty \int_0^1
\Tr \left(  
\frac{d\nabla(\epsilon)}{d \epsilon}
\: D \:
e^{- \: s^2 \: D^2} \right) \: du \: ds \:
 = \: - \: \frac{\sqrt{\pi}}{2} \:
d \:
\Tr\left(  
\frac{d\nabla(\epsilon)}{d \epsilon}
\: \frac{D}{|D|} \right).
\end{equation}

Using Proposition \ref{prop4} (in the finite-dimensional setting)
and the spectral theorem, we can deform $D$ to 
$P_+ \: - \: P_-$ without changing $\widetilde{\eta}^{(2)}$. Hence
we assume that $D \: = \: P_+ \: - \: P_-$.
Let us write $\nabla \: = \: \nabla_1 \: + \: \nabla_2$, where
$\nabla_1$ commutes with $D$ and $\nabla_2$ anticommutes with
$D$. Put $\nabla(\epsilon) \: = \: 
\nabla_1 \: + \: \epsilon \: \nabla_2$. Then
\begin{equation} \label{eqn4.10}
\Tr\left(  
\frac{d\nabla(\epsilon)}{d \epsilon}
\: \frac{D}{|D|} \right) \: = \: 
\Tr\left(  \nabla_2 \: 
D \right) \: = \: 0.
\end{equation}
Hence from (\ref{eqn4.9}), it suffices to compute
$\widetilde{\eta}^{(2)}$ when
$\nabla \: = \: \nabla_1$. In this case,
\begin{equation} \label{eqn4.11}
\widetilde{\eta}
\: = \: \int_0^\infty
\Tr_\sigma \left( \sigma \: D 
 \: 
e^{- \: s^2 \: D^2 \: - \: \nabla_1^2} \right) \: ds \: = \:
\frac{\sqrt{\pi}}{2} \: \Tr \left( \frac{D}{|D|} \: e^{- \: \nabla_1^2} 
\right),
\end{equation}
from which the proposition follows.
\end{pf}

\subsection{Infinite-dimensional case}

\begin{proposition} \label{prop6}
For $\epsilon \in [0,1]$, put 
\begin{equation} \label{eqn4.12}
A_s(\epsilon) \: = \: s \: \sigma \: D_\alpha 
\: + \: \nabla^{\pi_* E} \: + \: \frac{\epsilon}{4s} \: \sigma \: c(T).
\end{equation}
Define the corresponding $\widetilde{\eta}_\alpha(\epsilon) \in 
\Omega^{even}(B)$ as in (\ref{eqn2.19}).
Then given $D_0$, 
$\frac{d \widetilde{\eta}_\alpha^{(2)}(\epsilon)}{d\epsilon} \in
\Omega^2(B)$ is independent of the particular choice of $D_\alpha$.
\end{proposition}
\begin{pf}
We have $\frac{dA_s(\epsilon)}{d\epsilon} \: = \:
\frac{1}{4s} \: \sigma \: c(T)$.
From (\ref{eqn2.22}),
\begin{equation} \label{eqn4.13}
\frac{d \widetilde{\eta}_\alpha^{(2)}(\epsilon)}{d\epsilon} \: = \:
- \: \LIM_{s \rightarrow 0} \: 
\frac{1}{4s} \: \Tr \left( 
c(T) \:
e^{- s^2 \: D_\alpha^2} \right).
\end{equation}
At this point we do not have to
assume that $D_\alpha$ is invertible. Furthermore, the question is local
on $B$, so we may assume that $B \: = \: U_\alpha$.

If $D_\alpha \: = \: D_0$ then from standard heat equation asymptotics
\cite[Theorem 1.5]{Fegan-Gilkey (1985)},
\begin{equation} \label{eqn4.14}
\Tr \left( 
c(T) \:
e^{- s^2 \: D_0^2} \right) \: \sim \: s^{- \: dim(Z)} \: \sum_{k=0}^\infty
r_k \: s^{2k}
\end{equation}
for some $2$-form-valued coefficients $\{r_k\}_{k=0}^\infty$. Then
\begin{equation} \label{eqn4.15}
\frac{d \widetilde{\eta}_\alpha^{(2)}(\epsilon)}{d\epsilon} \: = \:
- \: \frac{1}{4} \: r_{\frac{dim(Z) \: + \: 1}{2}}.
\end{equation}
In particular, this may be nonzero.

For general $D_\alpha$, a Duhamel expansion around
$e^{- s^2 \: D_0^2}$ as in 
the proof of Proposition \ref{prop1} 
shows that
\begin{equation} \label{eqn4.16}
\LIM_{s \rightarrow 0} \: 
\frac{1}{4s} \: \Tr \left( 
c(T) \:
e^{- s^2 \: D_\alpha^2} \right) \: = \:
\LIM_{s \rightarrow 0} \: 
\frac{1}{4s} \: \Tr \left( 
c(T) \:
e^{- s^2 \: D_0^2} \right),
\end{equation}
from which the proposition follows.
\end{pf}

We now make the assumption that if $U_\alpha \cap U_\beta \: \neq \:
\emptyset$ then $D_\alpha$ commutes with $D_\beta$.
We wish to compute
$\widetilde{\eta}_\beta^{(2)} \: - \:
\widetilde{\eta}_\alpha^{(2)} \in \Omega^2(U_\alpha \cap U_\beta)$. 

Define 
$E_{-+}$ and $E_{+-}$ as in the introduction.
As $D_\beta \: - \: D_\alpha$ is smoothing,
$E_{-+}$ and $E_{+-}$ are finite-rank operators.
\begin{proposition} \label{prop7}
\begin{equation} \label{eqn4.17}
\widetilde{\eta}_\beta^{(2)} \: - \:
\widetilde{\eta}_\alpha^{(2)} \: = \: - \: \sqrt{\pi} \:
\Tr \left( \left( E_{-+} \: \nabla \: E_{-+} \right)^2 \: - \:
\left( E_{+-} \: \nabla \: E_{+-} \right)^2 \right).
\end{equation}
\end{proposition}
\begin{pf}
From Proposition \ref{prop6}, we may assume that the superconnection on 
$U_\alpha$ is
$s \: \sigma \: D_\alpha \: + \: \nabla$, and similarly on
$U_\beta$. 
For $\epsilon \in [0,1]$, put
\begin{align} \label{eqn4.18}
\nabla(\epsilon) \: = \: \nabla \: - \: \epsilon \: [ & 
(1 \: - \: E_{+-} \: - \: E_{-+} ) \: \nabla \: E_{+-} \: + \:
 (1 \: - \: E_{+-} \: - \: E_{-+} ) \: \nabla \: E_{-+} \: + \\
  & E_{+-} \: \: \nabla \:  (1 \: - \: E_{+-} \: - \: E_{-+} ) \: + \:
  E_{+-} \: \: \nabla \:  E_{-+} \: +  \notag \\
  &  E_{-+} \: \: \nabla \:  (1 \: - \: E_{+-} \: - \: E_{-+} ) \: + \:
  E_{-+} \: \: \nabla \:  E_{+-} ]. \notag
\end{align}
Then 
\begin{equation} \label{eqn4.19}
\nabla(1) \: = \: (1 \: - \: E_{+-} \: - \: E_{-+} ) \: \nabla \: 
(1 \: - \: E_{+-} \: - \: E_{-+} ) \: + \:
E_{+-} \: \nabla \: E_{+-} \: + \:
E_{-+} \: \nabla \: E_{-+}.
\end{equation} 

Note that $\frac{d\nabla(\epsilon)}{d\epsilon}$ is smoothing.
From (\ref{eqn2.22}),
\begin{align} \label{eqn4.20}
\frac{d \widetilde{\eta}_\alpha(\epsilon)}{d\epsilon}  \: = \:
& - \: \LIM_{s \rightarrow 0} \: 
\Tr_\sigma \left( 
\frac{d\nabla(\epsilon)}{d \epsilon} \: 
e^{- \: A_s(\epsilon)^2} \right) \: + \: \\
& d \: \LIM_{t \rightarrow 0} \int_t^\infty \int_0^1
\Tr_\sigma \left(  
\frac{d \nabla(\epsilon)}{d \epsilon}
\: e^{- \: u \: A_s(\epsilon)^2} \: 
\sigma \: D_\alpha \:
e^{- \: (1-u) \: A_s(\epsilon)^2} \right) \: du \: ds, \notag
\end{align}
and similarly for $\frac{d \widetilde{\eta}_\beta(\epsilon)}{d\epsilon}$.
Using the method of proof of Proposition \ref{prop1}, one finds
\begin{equation} \label{eqn4.21}
\LIM_{s \rightarrow 0} \: 
\Tr_\sigma \left( 
\frac{d\nabla(\epsilon)}{d \epsilon} \: 
e^{- \: A_s(\epsilon)^2} \right) \: = \: 0.
\end{equation}
Then as in (\ref{eqn4.9}),
\begin{align} \label{eqn4.22}
\frac{d \widetilde{\eta}_\alpha(\epsilon)^{(2)}}{d\epsilon}  \:  = \: - \:
d \: \int_0^\infty \int_0^1
\Tr \left(  
\frac{d\nabla(\epsilon)}{d \epsilon}
\: D_\alpha \:
e^{- \: s^2 \: D_\alpha^2} \right) \: du \: ds.
\end{align}
Thus
\begin{align} \label{eqn4.23}
\frac{d (\widetilde{\eta}_\beta(\epsilon)^{(2)} \: - \:
\widetilde{\eta}_\alpha(\epsilon)^{(2)})}{d\epsilon}
& = \: - \:
d \: \int_0^\infty
\Tr \left(  
\frac{d\nabla(\epsilon)}{d \epsilon}
\: \left( D_\beta \:
e^{- \: s^2 \: D_\beta^2} \: - \:
D_\alpha \:
e^{- \: s^2 \: D_\alpha^2}  
\right)
\right) \: ds \\
& = \: - \: \frac{\sqrt{\pi}}{2} \:
d \:
\Tr\left(  
\frac{d\nabla(\epsilon)}{d \epsilon}
\: \left( \frac{D_\beta}{|D_\beta|} \: - \:  
\frac{D_\alpha}{|D_\alpha|} \right) \right)  \notag \\
& = \: - \: \frac{\sqrt{\pi}}{2} \:
d \:
\Tr\left(  
\frac{d\nabla(\epsilon)}{d \epsilon}
\: \left( 2 \: E_{-+} \: - \: 2 \: E_{+-} \right) \right) \: = \: 0, \notag
\end{align}
where the last line comes from the off-diagonal nature of
$\frac{d\nabla(\epsilon)}{d \epsilon}$.
Hence we may assume that the superconnection on $U_\alpha$ is
$s \: \sigma \: D_\alpha \: + \: \nabla(1)$ and that
the superconnection on $U_\beta$ is
$s \: \sigma \: D_\beta \: + \: \nabla(1)$.

Due to the diagonal form (\ref{eqn4.19}) of $\nabla(1)$, 
$\widetilde{\eta}_\beta \: - \: \widetilde{\eta}_\alpha$ is the
sum of contributions from 
$\Image(I \: - \: E_{+-} \: - \: E_{-+})$ and
$\Image(E_{+-} \: + \: E_{-+})$.
From the method of proof of Proposition \ref{prop4}, the
contribution from $\Image(I \: - \: E_{+-} \: - \: E_{-+})$ vanishes.
Then
$\widetilde{\eta}_\beta \: - \: \widetilde{\eta}_\alpha$ is 
the difference
of the two eta-forms of the finite-dimensional vector bundle
$\Image(E_{-+}) \oplus \Image(E_{+-})$, equipped with the 
connection $E_{-+} \: \nabla \: E_{-+} \: \oplus \:
E_{+-} \: \nabla \: E_{+-}$, where the eta-form is computed first with
$D_\beta$ and then with $D_\alpha$.
The proposition now follows from
Proposition \ref{prop5}.
\end{pf}

Define $L_{\alpha \beta}$ as in (\ref{eqn1.1}), with its 
connection $\nabla_{\alpha \beta}$ induced from the
connections $E_{-+} \: \nabla \: E_{-+}$ and
$E_{+-} \: \nabla \: E_{+-}$. Let $F_{\alpha \beta}$ denote its curvature,
an imaginary-valued $2$-form on $U_\alpha \cap U_\beta$.
Proposition \ref{prop7} says that on 
$U_\alpha \cap U_\beta$,
\begin{equation} \label{eqn4.24}
\widehat{\eta}_\beta^{(2)} \: - \:
\widehat{\eta}_\alpha^{(2)} \: = \: - \:
\frac{F_{\alpha \beta}}{2\pi i}.
\end{equation}

Suppose that $U_\alpha \cap U_\beta \cap U_\gamma \: \neq \: \emptyset$.
Given $a, b, c \in \{+ 1, - 1\}$
let $H_{abc}$ be the subbundle of $\pi_* E$ on which
$\frac{D_\alpha}{|D_\alpha|}$ acts as multiplication by $a$,
$\frac{D_\beta}{|D_\beta|}$ acts as multiplication by $b$ and
$\frac{D_\gamma}{|D_\gamma|}$ acts as multiplication by $c$. Then
\begin{align} \label{eqn4.25}
L_{\alpha \beta} \: & \cong \: \Lambda^{max}(H_{-++}) \: \otimes \: 
\Lambda^{max}(H_{-+-}) \: \otimes \:
\left( \Lambda^{max}(H_{+-+}) \right)^{-1} \: \otimes \: 
\left( \Lambda^{max}(H_{+--}) \right)^{-1}, \\
L_{\beta \gamma} \: & \cong \: \Lambda^{max}(H_{+-+}) \: \otimes \: 
\Lambda^{max}(H_{--+}) \: \otimes \:
\left( \Lambda^{max}(H_{++-}) \right)^{-1} \: \otimes \: 
\left( \Lambda^{max}(H_{-+-}) \right)^{-1}, 
\notag \\
L_{\gamma \alpha} \: & \cong \: \Lambda^{max}(H_{++-}) \: \otimes \: 
\Lambda^{max}(H_{+--}) \: \otimes \:
\left( \Lambda^{max}(H_{-++}) \right)^{-1} \: \otimes \: 
\left( \Lambda^{max}(H_{--+}) \right)^{-1}. \notag
\end{align}
There is an obvious nowhere-zero section $\theta_{\alpha \beta \gamma}$ of 
$L_{\alpha \beta} \: \otimes \: L_{\beta \gamma} \: \otimes \: 
L_{\gamma \alpha}$ on $U_\alpha \cap U_\beta \cap U_\gamma$. 

In general, let $E$ be a Hermitian vector bundle with Hermitian
connection $\nabla^E$ and
let $E_1$ and $E_2$ be subbundles of $E$ such that there is an
orthogonal direct sum $E\: = \: E_1 \: \oplus \:
E_2$. We do not assume that $\nabla^E$ is diagonal with respect to $E_1$ and
$E_2$. Let $\nabla^{E_1}$ and $\nabla^{E_2}$ be the induced connections
on $E_1$ and $E_2$, respectively.  We have corresponding connections
$\nabla^{\Lambda^{max}(E)}$, $\nabla^{\Lambda^{max}(E_1)}$ and
$\nabla^{\Lambda^{max}(E_2)}$, on $\Lambda^{max}(E)$,
$\Lambda^{max}(E_1)$ and $\Lambda^{max}(E_2)$, respectively. 
Then with respect to the isomorphism
$\Lambda^{max}(E) \: \cong \: \Lambda^{max}(E_1) \: \otimes \: 
\Lambda^{max}(E_2)$, one can check that 
$\nabla^{\Lambda^{max}(E)} \: \cong \: \nabla^{\Lambda^{max}(E_1)} \: \otimes
\: \nabla^{\Lambda^{max}(E_2)}$. 
Recalling the definition of $L_{\alpha \beta}$ from (\ref{eqn1.1}),
it follows that the section 
$\theta_{\alpha \beta \gamma}$ is covariantly-constant with respect to
$\nabla_{\alpha \beta} \: \otimes \: \nabla_{\beta \gamma} \: \otimes \:
\nabla_{\gamma \alpha}$. 

Finally, if $U_\alpha \cap U_\beta \cap U_\gamma \cap
U_\delta \: \neq \: \emptyset$ then the cocycle condition
$\theta_{\beta \gamma \delta} \:
\theta_{\alpha \gamma \delta}^{-1} \:
\theta_{\alpha \beta \delta} \: 
\theta_{\alpha \beta \gamma}^{-1} \: = \: 1$ is obviously satisfied.

In summary, we have shown that $\{D_\alpha\}_{\alpha \in I}$ determine
a gerbe on $B$ with connection. From (\ref{eqn2.26}), its curvature is
\begin{equation} \label{eqn4.26}
\left( \int_Z \widehat{A} \left(R^{TZ}/2\pi i \right)
\: \wedge \: \ch \left( F^V/2\pi i \right) \right)^{(3)} \in \Omega^3(B).
\end{equation}
To recall, the conditions 
that we imposed on $\{D_\alpha\}_{\alpha \in I}$
were \\
1. $D_\alpha$ is invertible.\\
2. $D_\alpha \: - \: D_0$ is a smoothing operator.\\
3. If $U_\alpha \cap U_\beta \: \neq \: \emptyset$ then
$D_\alpha$ and $D_\beta$ commute. 

In order to construct a gerbe-with-connection that only depends on
$D_0$, $\nabla^{\pi_*E}$ and $T$, as in the introduction
we assume that $D_\alpha \: = \: D_0 \: + \: 
h_\alpha(D_0)$ for
some $h_\alpha \in C_c^\infty(\R)$.
From spectral theory and the continuity of the spectral projections 
of $(D_0)_b$ with
respect to $b \in B$, it is easy to see that such
$\{U_\alpha\}_{\alpha \in I}$ and $\{h_\alpha\}_{\alpha \in I}$ exist.
Given $\{U_\alpha\}_{\alpha \in I}$, 
suppose that we make another choice $\{h^\prime_\alpha\}_{\alpha \in I}$.
Put $D^\prime_\alpha \: = \: D_0 \: + \: h^\prime_\alpha(D_0)$. Note that
$D^\prime_\alpha$ commutes with $D_\alpha$. Let
$L^\prime_{\alpha \beta}$ denote the new line bundle with
connection on $U_\alpha \cap U_\beta$. Let $\widehat{\eta}^\prime_\alpha$
denote the new eta-form.

Define a line bundle $L_{\alpha \alpha^\prime}$, with connection, on 
$U_\alpha$ as in
(\ref{eqn1.1}), replacing the pair $(D_\alpha, D_\beta)$ by the
pair $(D_\alpha, D^\prime_\alpha)$. Let $F_{\alpha \alpha^\prime}$ denote
the corresponding curvature.  Then it follows that
\begin{equation} \label{eqn4.27}
L^\prime_{\alpha \beta} \: = \: (L_{\alpha \alpha^\prime})^{-1} \: \otimes \:
L_{\alpha \beta} \: \otimes \:  L_{\beta \beta^\prime}
\end{equation}
as a line bundle with connection. Furthermore, from Proposition \ref{prop7},
\begin{equation} \label{eqn4.28}
(\widehat{\eta}^\prime_\alpha)^{(2)} \: - \:
\widehat{\eta}_\alpha^{(2)} \: = \: - \: 
\frac{F_{\alpha \alpha^\prime}}{2\pi i}.
\end{equation}
Thus if we choose $\{h^\prime_\alpha\}_{\alpha \in I}$ instead of
$\{h_\alpha\}_{\alpha \in I}$, we obtain an equivalent connection.

Finally, if $\{U_\alpha\}_{\alpha \in I}$ and
$\{U^\prime_{\alpha^\prime}\}_{\alpha^\prime \in I^\prime}$ are two
open coverings then by taking a common refinement, we see that we obtain
isomorphic gerbes and connections. 
This proves Theorem \ref{thm1}.

It follows that (\ref{eqn4.26}), as a rational cohomology class, lies in the
image of $\HH^3(B; \Z) \rightarrow \HH^3(B; \Q)$. Of course, one can see
this directly.
\begin{proposition} \label{prop8}
(\ref{eqn4.26}), as a rational cohomology class, lies in the
image of $\HH^3(B; \Z) \rightarrow \HH^3(B; \Q)$.
\end{proposition}
\begin{pf}
From the universal coefficient theorem, it is enough to show that the result
of pairing the rational cohomology class
(\ref{eqn4.26}) with an integer homology class $y \in \HH_3(B; \Z)$
is an integer.
As the map $s \: : \: \Omega^{SO}_3(B) \rightarrow \HH_3(B; \Z)$ from oriented
bordism to integer homology is surjective, we may assume that there is
a closed oriented $3$-manifold $X$ and a smooth map 
$\phi \: : \: X \rightarrow B$ such that $y \: = \: \phi_*([X])$, where
$[X] \in \HH_3(X; \Z)$ is the fundamental class of $X$. Then we can compute
the pairing of (\ref{eqn4.26}) with $y$ by pulling back 
(\ref{eqn4.26}) under $\phi$ to
$X$ and computing its pairing with $[X]$. Let $\pi^\prime \: : \:
M^\prime \rightarrow X$
denote the fiber bundle obtained by pulling back the fiber bundle
$\pi \: : \: M \rightarrow B$
under $\phi \: : \: X \rightarrow B$. Let $Z^\prime$ denote the fiber of
$\pi^\prime \: : \:
M^\prime \rightarrow X$ and let $V^\prime$ denote the pullback of $V$ to
$M^\prime$. Then by naturality, it is enough to
show that 
\begin{equation} \label{eqn4.29}
\int_X \left( \int_{Z^\prime} \widehat{A} \left(R^{TZ^\prime}/2\pi i \right)
\: \wedge \: \ch \left( F^{V^\prime}/2\pi i \right) \right)
\end{equation}
is an integer.

As $TZ$ has a spin structure, $TZ^\prime$ has a spin structure.
As $X$ is an oriented $3$-manifold, it has a spin structure.
Then $TM^\prime \: = \: TZ^\prime \: \oplus \: (\pi^\prime)^* TX$ has a
spin structure. Furthermore
\begin{equation} \label{eqn4.30}
\widehat{A} \left(R^{TM^\prime}/2\pi i \right) \: = \: 
\widehat{A} \left(R^{TZ^\prime}/2\pi i \right) \: \cup \: 
(\pi^\prime)^* \widehat{A} \left(R^{TX}/2\pi i \right) \: = \: 
\widehat{A} \left(R^{TZ^\prime}/2\pi i \right).
\end{equation}
Hence
\begin{equation} \label{eqn4.31}
\int_X \left( \int_{Z^\prime} \widehat{A} \left(R^{TZ^\prime}/2\pi i \right)
\: \wedge \: \ch \left( F^{V^\prime}/2\pi i \right) \right) \: = \:
\int_{M^\prime} \widehat{A} \left(R^{TM^\prime}/2\pi i \right)
\: \wedge \: \ch \left( F^{V^\prime}/2\pi i \right).
\end{equation}
The right-hand-side of (\ref{eqn4.31}) is an integer by the Atiyah-Singer index
theorem.
\end{pf}

\section{Deligne Cocycles} \label{Deligne}

We now assume that $Z$, the fiber of the fiber bundle $\pi \: : \: M
\rightarrow B$, is even-dimensional. We let $D$ denote the ensuing
family $D \: = \: \{D_b\}_{b \in B}$ of Dirac-type operators, with 
$D_b$ acting on
$C^\infty \left(Z_b; E \big|_{Z_b} \right)$. (We previously called this
$D_0$.)

Given $s \: > \: 0$, the corresponding Bismut superconnection 
\cite[Section III]{Bismut (1986)}, 
\cite[Chapter 10.3]{Berline-Getzler-Vergne (1992)}
on 
$\pi_* E$ is
\begin{equation} \label{eqn5.1}
A_s \: = \: s \: D
\: + \: \nabla^{\pi_* E} \: + \: \frac{1}{4s} \: c(T).
\end{equation}
For any $s \: > \: 0$, the supertrace
$\Tr_s \left( e^{- \: A_s^2} \right) \in \Omega^{even}(B)$
represents the Chern character of the index $\Ind(D) \in \KK^0(B)$ of
the family of vertical operators, up
to normalizing constants.

We have 
\cite[Theorem 10.23]{Berline-Getzler-Vergne (1992)},\cite{Bismut (1986)}
\begin{equation} \label{eqn5.2}
\lim_{s \rightarrow 0} \Tr_s \left( e^{- \: A_{0,s}^2} \right) \: = \:
(2 \pi i)^{- \:
\frac{dim(Z)}{2}} \:
\int_Z \widehat{A} \left(R^{TZ} \right)
\: \wedge \: \ch \left( F^V \right).
\end{equation}

\subsection{The case of vector bundle kernel}

We now make the assumption that $\Ker(D)$ has constant rank, i.e. is a
$\Z_2$-graded vector bundle on $B$. 
We give $\Ker(D)$ the projected connection
$\nabla^{Ker(D)}$ from $\nabla^{\pi_*E}$. It preserves the $\Z_2$-grading
on $\Ker(D)$. Let $F^{Ker(D)}$ denote the curvature of
$\nabla^{Ker(D)}$.

As in
\cite[Definition 4.33]{Bismut-Cheeger (1989)} and 
\cite[p. 273]{Dai (1991)}, define
$\widetilde{\eta} \in \Omega^{odd}(B)$ by
\begin{equation} \label{eqn5.3}
\widetilde{\eta}
\: = \: \int_0^\infty
\Tr_s \left( 
\frac{d A_s}{d s} \: 
e^{- \: A_s^2} \right) \: ds.
\end{equation}
Then from \cite[Theorem 4.35]{Bismut-Cheeger (1989)} and 
\cite[Theorem 0.1]{Dai (1991)},
\begin{equation} \label{eqn5.4}
d \widetilde{\eta} \: = \: (2 \pi i)^{- \:
\frac{dim(Z)}{2}} \:
\int_Z \widehat{A} \left( R^{TZ} \right)
\: \wedge \: \ch \left( F^V \right) \: - \: \ch \left( F^{Ker(D)} \right).
\end{equation}

To give normalizations that are compatible with rational cohomology,
let ${\cal R}$ be the operator on $\Omega^*(B)$ which acts on
$\Omega^{2k}(B)$ as multiplication by $(2\pi i)^{-k}$ and which acts
on $\Omega^{2k+1}(B)$ as multiplication by $(2\pi i)^{-k}$.
Put 
\begin{equation} \label{eqn5.5}
\widehat{\eta} \: = \: \frac{1}{2 \pi i} \: {\cal R} \: \widetilde{\eta}.
\end{equation}
Then (\ref{eqn5.4}) becomes
\begin{equation} \label{eqn5.6} 
d \widehat{\eta} \: = \:
\int_Z \widehat{A} \left(R^{TZ}/2\pi i \right)
\: \wedge \: \ch \left( F^V/2\pi i \right)
 \: - \: \ch \left( F^{Ker(D)}/ 2\pi i \right).
\end{equation}

We now make the assumption that for each $b \in B$, the index of
$D_b$ vanishes in $\Z$.  Equivalently, the vector bundles $\Ker(D)_+$ and
$\Ker(D)_-$ have the same rank. For simplicity of notation, we will
abbreviate $\nabla^{Ker(D)_\pm}$ by $\nabla^\pm$, and write its curvature as
$F^\pm$.

Let $\{U_\alpha\}_{\alpha \in I}$ 
be a covering of $B$ by open sets such that over $U_\alpha$,
there is an isometric isomorphism $W_\alpha \: : \: 
\Ker(D)_+ \big|_{U_\alpha} \rightarrow \Ker(D)_- \big|_{U_\alpha}$.
For example, if each $U_\alpha$ is contractible then such $W_\alpha$'s exist.

Define the Chern-Simons form $\CS_\alpha \in \Omega^{odd}(U_\alpha)$ by
\begin{equation} \label{eqn5.7}
\CS_\alpha \: = \: - \: \frac{1}{2 \pi i} \:
\int_0^1 \Tr \left( \nabla^{+} \: - \: W_\alpha^{-1} \circ
\nabla^{-} \circ W_\alpha \right) \: 
e^{- \: \frac{\left(t \nabla^{+} \: + \: (1 \: - \: t ) \:
W_\alpha^{-1} \circ \nabla^{-} \circ W_\alpha  \right)^2}{2 \pi i}} \: dt.
\end{equation}
By construction,
\begin{equation} \label{eqn5.8}
d \CS_\alpha \: = \: \ch \left( F^+ / 2 \pi i \right) \: - \: 
\ch \left( W_\alpha^{-1} \circ F^- \circ W_\alpha / 2 \pi i \right)
\: = \: 
\ch \left( F^{Ker(D)} / 2 \pi i \right).
\end{equation}
Then from (\ref{eqn5.6}),
\begin{equation} \label{eqn5.9}
d \left( \widehat{\eta} \: + \: \CS_\alpha \right) \: = \: 
\int_Z \widehat{A} \left(R^{TZ}/2\pi i \right)
\: \wedge \: \ch \left( F^V/2\pi i \right).
\end{equation}
on $U_\alpha$.
Thus we wish to attach the odd form $\widehat{\eta} \: + \: \CS_\alpha$
to $U_\alpha$.

If $U_\alpha \cap U_\beta \neq \emptyset$ then we wish to write
$\left( \widehat{\eta} \: + \: \CS_\beta \right) \: - \: 
\left( \widehat{\eta} \: + \: \CS_\alpha \right) \: = \: 
\CS_\beta \: - \: \CS_\alpha$ as an exact form on $U_\alpha \cap U_\beta$.
Then we wish to repeat the process if $U_\alpha \cap U_\beta \cap U_\gamma
\neq \emptyset$, etc.  In order to streamline things, we use a 
construction which is similar to the ``descent equations'' in the
study of anomalies \cite{Atiyah-Singer (1984),Zumino-Wu-Zee (1984)}.

For simplicity, we assume for the moment that $U_\alpha \: = \: B$.
We write ${\cal I} \: = \: \Isom(\Ker(D)_+, \Ker(D)_-)$. It is acted upon
freely and transitively by the groups of gauge transformations 
$\Isom(\Ker(D)_+)$ and $\Isom(\Ker(D)_-)$. We let $W$ denote a ``coordinate''
on ${\cal I}$ and we let $\delta$ denote the
differential on ${\cal I}$, so that
$W^{-1} \: \delta W$ denotes the canonical left-$\Isom(\Ker(D)_-)$-invariant
$1$-form on ${\cal I}$, with values in
$\End(\Ker(D)_+)$. 

Consider the vector bundle ${\cal I} \: \times \: B \: 
\times \: \Ker(D)_+$ on ${\cal I} \: \times \: B$.
It has two canonical connections, $\delta \: + \: \nabla^+$ and
$W^{-1} \circ (\delta \: + \: \nabla^-) \circ W \: = \:
\delta \: + \: W^{-1} \: \delta W \: + \: 
W^{-1} \circ \nabla^- \circ W$.
Consider the corresponding Chern-Simons form
$\CS \in \Omega^{odd}({\cal I} \: \times \: B)$
given by
\begin{align} \label{eqn5.10}
\CS \: & = \: - \: \frac{1}{2 \pi i} \:
\int_0^1 \Tr \left( 
\delta \: + \: \nabla^{+} \: - \: W^{-1} \circ (\delta \: + \:
\nabla^{-}) \circ W \right) \: 
e^{- \: \frac{\left(t( \delta \: + \: \nabla^{+}) \: + \: (1 \: - \: t ) \:
W^{-1} \circ (\delta \: + \: 
\nabla^{-}) \circ W \right)^2}{2 \pi i}} \: dt \\
& = \: - \: \frac{1}{2 \pi i} \:
\int_0^1 \Tr \left( \nabla^{+} \: - \: W^{-1} \: \delta W \: - \:
W^{-1} \circ
\nabla^{-} \circ W \right) \: 
e^{- \: \frac{\left(\delta \: + \: t \nabla^{+} \: + \: (1 \: - \: t ) \:
(W^{-1} \: \delta W \: + \: W^{-1} \circ 
\nabla^{-} \circ W ) \right)^2}{2 \pi i}} \: dt \notag.
\end{align}
By construction,
\begin{equation} \label{eqn5.11}
(\delta \: + \: d) \: \CS \: = \: 
\ch \left( F^{+} / 2 \pi i \right) \: - \: 
\ch \left( W^{-1} \circ F^{-} \circ W / 2 \pi i \right)
\: = \: 
\ch \left( F^{Ker(D)} / 2 \pi i \right).
\end{equation}

Let us work out $\CS$ in low degrees. To do so, we use the fact that
\begin{align} \label{eqn5.12}
& \left(\delta \: + \: t \nabla^{+} \: + \: (1 \: - \: t ) \:
(W^{-1} \: \delta W \: + \: W^{-1} \circ 
\nabla^{-} \circ W \right)^2 \: = \\
& t \: F^+ \: + \: (1 \: - \: t) \: W^{-1} \circ F^- \circ W \: - \:
t \: (1 \: - \: t) \: 
\left( \nabla^{+} \: - \: W^{-1} \: \delta W \: - \:
W^{-1} \circ
\nabla^{-} \circ W \right)^2. \notag
\end{align}

\subsubsection{$\HH^2$}

The $1$-form component of $\CS$ is
\begin{equation} \label{eqn5.13}
\CS^{(1)} \: = \: - \: \frac{1}{2 \pi i} \:
\Tr \left( 
\nabla^{+} \: - \: W^{-1} \: \delta W \: - \:
W^{-1} \circ \nabla^{-} \circ W \right).
\end{equation}
Then (\ref{eqn5.11}) becomes
\begin{align} \label{eqn5.14}
d \left( - \: \frac{1}{2 \pi i} \:
\Tr ( \nabla^+ \: - \: W^{-1} \circ \nabla^- \circ W)
  \right) \: & = \: c_1(\nabla^+) \: - \: c_1(\nabla^-), \\
- \: d \: \Tr (W^{-1} \: \delta W) \: + \: \delta \:
\Tr \left( \nabla^{+} \: - \:
W^{-1} \circ
\nabla^{-} \circ W \right) \: & = \: 0,  \notag \\
\delta \: \Tr \left( W^{-1} \:\delta W \right) \: & = \: 0. \notag
\end{align}

Given $\{U_\alpha\}_{\alpha \in I}$ and $\{W_\alpha\}_{\alpha \in I}$ as
before, suppose that $U_\alpha \cap U_\beta \neq \emptyset$.
Let ${\cal I}_{\alpha \beta}$ be the space ${\cal I}$ defined above when
the base is $U_\alpha \cap U_\beta$. Suppose that
there is a smooth path $\sigma_{\alpha \beta} \: : \: [0,1] \rightarrow 
{\cal I}_{\alpha \beta}$ from 
$W_\alpha \big|_{U_\alpha \cap U_\beta}$ to 
$W_\beta \big|_{U_\alpha \cap U_\beta}$. For example,
if $U_\alpha \cap U_\beta$ is contractible then there is such a path,
as the unitary group is connected. We put 
$\sigma_{\beta \alpha}(t) \: = \: 
\sigma_{\alpha \beta}(1 \: - \: t)$.
It makes sense to write
\begin{equation} \label{eqn5.15}
\int_{[0,1]} \sigma_{\alpha \beta}^* \:  \Tr (W^{-1} \: \delta W) \: = \:
\int_{[0,1]} \Tr \left( \sigma_{\alpha \beta}(t)^{-1} \: 
\frac{d \sigma_{\alpha \beta}(t)}{dt} \right) \: dt \in
\Omega^0(U_\alpha \cap U_\beta)
\end{equation}
and from (\ref{eqn5.14}),
\begin{align} \label{eqn5.16}
& d \left( \frac{1}{2 \pi i} \: \int_{[0,1]} \sigma_{\alpha \beta}^* \: 
\Tr \left( W^{-1} \: \delta W \right) \right) \: = \\
& 
- \: \frac{1}{2 \pi i} \: \Tr \: \left(
\nabla^{+} \: - \: W_\beta^{-1} \circ \nabla^{-} \circ W_\beta \right) \: + \: 
 \: \frac{1}{2 \pi i} \: \Tr \: \left(
\nabla^{+} \: - \: W_\alpha^{-1} \circ \nabla^{-} \circ W_\alpha \right) 
\: = \: \notag \\
& \CS_\beta^{(1)} \: - \: \CS_\alpha^{(1)}.
\end{align}

If $U_\alpha \cap U_\beta \cap U_\gamma \neq \emptyset$, let
${\cal I}_{\alpha \beta \gamma}$ be the space ${\cal I}$ defined as
above when the base is $U_\alpha \cap U_\beta \cap U_\gamma$. Let
$\mu_{\alpha \beta \gamma} \: : \: S^1 \rightarrow 
{\cal I}_{\alpha \beta \gamma}$ be a smooth concatenation of
$\sigma_{\alpha \beta}$, $\sigma_{\beta \gamma}$ and
$\sigma_{\gamma \alpha}$. Then
$\frac{1}{2 \pi i} \:
\int_{S^1} \mu_{\alpha \beta \gamma}^*\: \Tr (W^{-1} \: \delta W)$
is a continuous integer-valued function on
$U_\alpha \cap U_\beta \cap U_\gamma$.

We can summarize the discussion so far by saying that
\begin{equation} \label{eqn5.17}
{\cal C} \: = \: \left( \widehat{\eta}^{(1)} \: + \: \CS_\alpha^{(1)}, \:
\frac{1}{2 \pi i} \:
\int_{[0,1]} \sigma_{\alpha \beta}^*\: \Tr (W^{-1} \: \delta W), \:
\frac{1}{2 \pi i} \:
\int_{S^1} \mu_{\alpha \beta \gamma}^*\: \Tr (W^{-1} \: \delta W)
 \right)
\end{equation}
forms a 2-cocycle for the \v{C}ech-cohomology of the
complex of sheaves
\begin{equation} \label{eqn5.18}
\Z \longrightarrow \Omega^0 \longrightarrow \Omega^1
\end{equation} 
on $B$,
where $\Omega^p$ denotes the sheaf of real-valued $p$-forms.

Now suppose that $\{W^\prime_\alpha\}_{\alpha \in I}$ is another choice
of isometries, with each $W^\prime_\alpha$ connected to
$W_\alpha$ in ${\cal I}_\alpha$. Let $\sigma_{\alpha^\prime \beta^\prime}$ be a
path from $W^\prime_\alpha \big|_{U_\alpha \cap U_\beta}$ to 
$W^\prime_\beta \big|_{U_\alpha \cap U_\beta}$. We obtain a 
corresponding cocycle
\begin{equation} \label{eqn5.19}
{\cal C}^\prime \: = \:
\left( \widehat{\eta}^{(1)} \: + \: \CS^{\prime,(1)}_\alpha, \:
\frac{1}{2 \pi i} \:
\int_{[0,1]} \sigma_{\alpha^\prime \beta^\prime}^*\: 
\Tr (W^{-1} \: \delta W), \:
\frac{1}{2 \pi i} \:
\int_{S^1} \mu_{\alpha^\prime \beta^\prime \gamma^\prime}^*\: 
\Tr (W^{-1} \: \delta W)
 \right).
\end{equation}
In order to compare ${\cal C}$ and ${\cal C}^\prime$, 
for each $\alpha \in I$ choose a path
$\sigma_{\alpha \alpha^\prime} \: : \: [0, 1] \rightarrow {\cal I}_\alpha$
from $W_\alpha$ to $W_\alpha^\prime$. 
If $U_\alpha \cap U_\beta \neq \emptyset$ then
define 
$\mu_{\alpha^\prime \beta^\prime \beta \alpha} \: : \: S^1 \:
\rightarrow {\cal I}_{\alpha \beta}$
in the obvious way. Then
\begin{align} \label{eqn5.20}
& ( \widehat{\eta} \: + \: \CS^{\prime,(1)}_\alpha) \: - \:
( \widehat{\eta} \: + \: \CS^{(1)}_\alpha) \:  = \: d \left(
\frac{1}{2 \pi i} \:
\int_{[0,1]} \sigma_{\alpha \alpha^\prime}^* \: 
\Tr (W^{-1} \: \delta W) \right), \\
& \frac{1}{2 \pi i} \:
\int_{[0,1]} \sigma_{\alpha^\prime \beta^\prime}^*
\: \Tr (W^{-1} \: \delta W) 
 \: - \: 
\frac{1}{2 \pi i} \:
\int_{[0,1]} \sigma_{\alpha \beta}^*\: \Tr (W^{-1} \: \delta W) 
\: = \notag \\
&\: \frac{1}{2 \pi i} \:
\int_{[0,1]} \sigma_{\beta \beta^\prime}^*\: \Tr (W^{-1} \: \delta W) 
\: - \: 
\frac{1}{2 \pi i} \:
\int_{[0,1]} \sigma_{\alpha \alpha^\prime}^*\: \Tr (W^{-1} \: \delta W) 
 \notag \\
& + \: \frac{1}{2 \pi i} \:
\int_{S^1} \mu_{\alpha^\prime \beta^\prime \beta \alpha}^*
\: \Tr (W^{-1} \: \delta W), \notag \\
& \frac{1}{2 \pi i} \:
\int_{S^1} \mu_{\alpha^\prime \beta^\prime \gamma^\prime}^*\: 
\Tr (W^{-1} \: \delta W) 
\: - \: 
\frac{1}{2 \pi i} \:
\int_{S^1} \mu_{\alpha \beta \gamma}^*\: \Tr (W^{-1} \: \delta W) 
\:  = \notag \\
& \frac{1}{2 \pi i} \: 
\int_{S^1} \left( \mu_{\alpha^\prime \beta^\prime \beta \alpha}^* \: 
\Tr (W^{-1} \: \delta W) 
 \: + \: \mu_{\beta^\prime \gamma^\prime \gamma \beta}^* \: 
\Tr (W^{-1} \: \delta W) \: + \:
\mu_{\gamma^\prime \alpha^\prime \alpha \gamma}^* \: 
\Tr (W^{-1} \: \delta W) 
\right). \notag
\end{align}
In other words, ${\cal C}^\prime \: - \: {\cal C}$ is the coboundary of
the $1$-cochain
\begin{equation} \label{eqn5.21}
\left( \frac{1}{2 \pi i} \: 
\int_{[0,1]} \sigma_{\alpha \alpha^\prime}^* \: 
\Tr (W^{-1} \: \delta W), \:
\frac{1}{2 \pi i} \: 
\int_{S^1}  \mu_{\alpha^\prime \beta^\prime \beta \alpha}^* \: 
\Tr (W^{-1} \: \delta W) 
\right).
\end{equation}
Thus ${\cal C}$ and ${\cal C}^\prime$ are cohomologous.

In summary, our input data consisted of points 
$W_\alpha \in {\cal I}_{\alpha}$,
defined up to homotopy, with the property that if $U_\alpha \cap U_\beta
\neq \emptyset$ then there is a path in ${\cal I}_{\alpha \beta}$ from
$W_\alpha \big|_{U_\alpha \cap U_\beta}$ to 
$W_\beta \big|_{U_\alpha \cap U_\beta}$.
From this we obtained a 
Deligne cohomology class on $B$ of degree $2$
\cite[Chapter 1.5]{Brylinski (1993)}.
From (\ref{eqn5.9}), its ``curvature'' is the $2$-form
\begin{equation} \label{eqn5.22}
\left( \int_Z \widehat{A} \left(R^{TZ}/2\pi i \right)
\: \wedge \: \ch \left( F^V/2\pi i \right)
 \right)^{(2)} \in
\Omega^2(B).
\end{equation}

In view of the isomorphism between the $2$-dimensional Deligne
cohomology of $B$ and the isomorphism classes of line bundles with
connection on $B$ \cite[Theorem 2.2.12]{Brylinski (1993)}, 
we also obtain a line bundle with connection on $B$ which is,
of course, the determinant line bundle. To see this explicitly, let us first
note that $e^{\int_{[0,1]} \sigma_{\alpha \beta}^* \:  Tr (W^{-1} \: \delta W)}
\in \U(1)$ depends only on $W_\alpha$ and $W_\beta$, and not on
$\sigma_{\alpha \beta}$. To evaluate it,
formally
\begin{align} \label{eqn5.23}
\int_{[0,1]} \sigma_{\alpha \beta}^* \:  \Tr (W^{-1} \: \delta W) \: & = \:
\int_{[0,1]} \sigma_{\alpha \beta}^* \:  \delta \ln \det(W) \: = \:
\int_{[0,1]} d \: \sigma_{\alpha \beta}^* \: \ln \det(W) \\
& = \: \ln \det(W_\beta) \: - \: \ln \det(W_\alpha) \: = \:
\ln \det(W_\alpha^{-1} \: W_\beta). \notag
\end{align}
These equations make sense modulo $2 \pi i \Z$, to give
\begin{equation} \label{eqn5.24}
e^{\int_{[0,1]} \sigma_{\alpha \beta}^* \:  Tr (W^{-1} \: \delta W)}
\: = \: \det(W_\alpha^{-1} \: W_\beta).
\end{equation}

Consider the imaginary-valued $1$-form
$A_\alpha \: =  - \: 2 \pi i \: 
\left( \widehat{\eta} \: + \: \CS_\alpha^{(1)} \right)$ on $U_\alpha$.
Equations (\ref{eqn5.20}) and (\ref{eqn5.24}) show that the forms
$\{A_\alpha\}_{\alpha \in I}$ fit together to give a connection on the
line bundle whose transition functions are
$\phi_{\alpha \beta} \: = \: \det(W_\alpha^{-1} \: W_\beta)$.
This is the same as the determinant line bundle
\cite[Chapter 9.7]{Berline-Getzler-Vergne (1992)}, which in our case is
equal to 
$\Lambda^{max}(\Ker(D)_+)^{-1} \: \otimes \: \Lambda^{max}(\Ker(D)_-)$.
The connection that we have defined on the determinant line bundle is the same
as that defined in \cite[Chapter 9.7]{Berline-Getzler-Vergne (1992)}.
Its curvature is given by (\ref{eqn5.22}).
Of course the determinant line bundle can be defined without the
assumption that $\Ker(D)$ is a vector bundle on $B$, or our other assumptions.

\subsubsection{$\HH^4$}

The $3$-form component of $\CS$ is
\begin{align} \label{eqn5.25}
\CS^{(3)} \: = \: & - \: \frac{1}{(2 \pi i)^2} \: \left[
\frac{1}{2} \: \Tr \: \left(
(\nabla^+ \: - \: W^{-1} \: \nabla^- \: W) \: \wedge \: \left( 
F^+ \: + \: W^{-1} \: F^- \: W \: \right. \right. \right.\\
& \left. \left.
\: \: \: \: \: \: \: \: \: \: \: \: \: \: \: \: \: \: \: \: \: \: 
- \: \frac{1}{3}
(\nabla^+ \: - \: W^{-1} \: \nabla^- \: W)^2 \right) \right) \notag \\
& - \: \frac{1}{2} \: \Tr \: \left( W^{-1} \: \delta W \wedge (F^+ \: + \: 
W^{-1} \: F^- \: W \: - \: 
(\nabla^+ \: - \: W^{-1} \: \nabla^- \: W)^2 ) \right) \notag \\
&-  \: \frac{1}{2} \: \Tr \: \left( W^{-1} \: \delta W \: \wedge \:  
W^{-1} \: \delta W \: \wedge \:
(\nabla^+ \: - \: W^{-1} \: \nabla^- \: W) \right) \notag \\
& \left. + \: \frac{1}{6} \:  \Tr \: \left(
W^{-1} \: \delta W \right)^3 \right]. \notag
\end{align}
Then (\ref{eqn5.11}) becomes
\begin{align} \label{eqn5.26}
& d \: 
 \Tr \: \left(
(\nabla^+ \: - \: W^{-1} \: \nabla^- \: W) \: \wedge \: 
\left( F^+ \: + \: W^{-1} \: F^- \: W \:
- \: \frac{1}{3}
(\nabla^+ \: - \: W^{-1} \: \nabla^- \: W)^2 \right) \right) \: = \\
& \Tr \left( F^+ \right)^2 \: - \: 
\Tr \left( F^- \right)^2, \notag \\
& - \: \frac{1}{2} \: d \: 
\Tr \: \left( W^{-1} \: \delta W \wedge (F^+ \: + \: 
W^{-1} \: F^- \: W \: - \: 
(\nabla^+ \: - \: W^{-1} \: \nabla^- \: W)^2 ) \right) \: + \notag \\
& \frac{1}{2} \: \delta \: 
 \Tr \: \left(
(\nabla^+ \: - \: W^{-1} \: \nabla^- \: W) \: \wedge \: 
\left( F^+ \: + \: W^{-1} \: F^- \: W \:
- \: \frac{1}{3}
(\nabla^+ \: - \: W^{-1} \: \nabla^- \: W)^2 \right) \right) \: = \: 0, 
\notag \\
&-  \: \frac{1}{2} \: d \: \Tr \: \left( W^{-1} \: \delta W \: \wedge \:  
W^{-1} \: \delta W \: \wedge \:
(\nabla^+ \: - \: W^{-1} \: \nabla^- \: W) \right) \notag \\
& - \: \frac{1}{2} \: \delta \: 
\Tr \: \left( W^{-1} \: \delta W \wedge (F^+ \: + \: 
W^{-1} \: F^- \: W \: - \: 
(\nabla^+ \: - \: W^{-1} \: \nabla^- \: W)^2 ) \right)
\: = \: 0,  \notag \\
& \frac{1}{6} \:  d \: \Tr \: \left(
W^{-1} \: \delta W \right)^3 \: 
-  \: \frac{1}{2} \: \delta \: \Tr \: \left( W^{-1} \: \delta W \: \wedge \:  
W^{-1} \: \delta W \: \wedge \:
(\nabla^+ \: - \: W^{-1} \: \nabla^- \: W) \right) \: = \: 0, \notag \\
& \delta \: \Tr \left( W^{-1} \: \delta W \right)^3 \: = \: 0. \notag
\end{align}

Given $\{U_\alpha\}_{\alpha \in I}$ 
we suppose that there exist \\
1. Maps $\sigma_\alpha \: : \: D^0 \rightarrow {\cal I}_\alpha$, \\
2. Maps $\sigma_{\alpha \beta} \: : \: D^1 \rightarrow 
{\cal I}_{\alpha \beta}$ such that 
\begin{equation} \label{eqn5.27}
\partial \sigma_{\alpha \beta} \: = \:
\sigma_\beta \big|_{U_\alpha \cap U_\beta} \: - \: 
\sigma_\alpha \big|_{U_\alpha \cap U_\beta},
\end{equation}
3. Maps $\sigma_{\alpha \beta \gamma} \: : \: D^2 \rightarrow
{\cal I}_{\alpha \beta \gamma}$ such that 
\begin{equation} \label{eqn5.28}
\partial \sigma_{\alpha \beta \gamma} \: = \:
\sigma_{\beta \gamma} \big|_{U_\alpha \cap U_\beta \cap U_\gamma} \: - \:
\sigma_{\alpha \gamma} \big|_{U_\alpha \cap U_\beta \cap U_\gamma} \: + \: 
\sigma_{\alpha \beta} \big|_{U_\alpha \cap U_\beta \cap U_\gamma}
\end{equation}
and \\
4. Maps $\sigma_{\alpha \beta \gamma \delta} \: : \: D^3 \rightarrow
{\cal I}_{\alpha \beta \gamma \delta}$ such that 
\begin{align} \label{eqn5.29}
\partial \sigma_{\alpha \beta \gamma \delta} \: = \:
& \sigma_{\beta \gamma \delta}
\big|_{U_\alpha \cap U_\beta \cap U_\gamma \cap U_\delta} \: 
- \: \sigma_{\alpha \gamma \delta} 
\big|_{U_\alpha \cap U_\beta \cap U_\gamma \cap U_\delta} \: + \\
& \sigma_{\alpha \beta \delta} 
\big|_{U_\alpha \cap U_\beta \cap U_\gamma \cap U_\delta} \: - \: 
\sigma_{\alpha \beta \gamma}
\big|_{U_\alpha \cap U_\beta \cap U_\gamma \cap U_\delta}. \notag
\end{align}

The right-hand-side of (\ref{eqn5.28}), for example, means a concatenation
$\sigma_{\alpha \beta} \cup \sigma_{\beta \gamma} \cup
\sigma_{\gamma \alpha}$.  {\it A priori} this is a continuous map from
$S^1$ to ${\cal I}_{\alpha \beta}$, but after an appropriate reparametrization
we may assume that it is smooth.
Equation (\ref{eqn5.28}) means that
$\sigma_{\alpha \beta} \cup \sigma_{\beta \gamma} \cup
\sigma_{\gamma \alpha}$ extends to a map from $D^2$ to
${\cal I}_{\alpha \beta}$. Again we may assume first that it extends
continuously and then obtain a smooth extension. (We could also work with
piecewise smooth maps.)

If $U_\alpha \cap U_\beta \cap U_\gamma \cap U_\delta \cap U_\epsilon \neq
\emptyset$, define
$\mu_{\alpha \beta \gamma \delta \epsilon} \: : S^3 \rightarrow 
{\cal I}_{\alpha \beta \gamma \delta \epsilon}$ by
\begin{align} \label{eqn5.30}
\mu_{\alpha \beta \gamma \delta \epsilon} \: = \: 
& \sigma_{\beta \gamma \delta \epsilon} 
\big|_{U_\alpha \cap U_\beta \cap U_\gamma \cap U_\delta 
\cap U_\epsilon} \: - \:
\sigma_{\alpha \gamma \delta \epsilon} 
\big|_{U_\alpha \cap U_\beta \cap U_\gamma \cap U_\delta
\cap U_\epsilon} \: + \:
 \sigma_{\alpha \beta \delta \epsilon}
\big|_{U_\alpha \cap U_\beta \cap U_\gamma \cap U_\delta
\cap U_\epsilon}  \: -  \\
& \sigma_{\alpha \beta \gamma \epsilon}
\big|_{U_\alpha \cap U_\beta \cap U_\gamma \cap U_\delta
\cap U_\epsilon}  \: + \:
\sigma_{\alpha \beta \gamma \delta}
\big|_{U_\alpha \cap U_\beta \cap U_\gamma \cap U_\delta
\cap U_\epsilon}. \notag
\end{align}

Put 
\begin{align} \label{eqn5.31}
G_\alpha \: = \: & \widehat{\eta}^{(3)} \: + \: \CS_\alpha^{(3)}, \\
G_{\alpha \beta} \: = \: & \frac{1}{2} \: \frac{1}{(2 \pi i)^2} \:
\int_{D^1} \sigma_{\alpha \beta}^* \: \Tr \left(
W^{-1} \: \delta W \: \wedge \: (F^+ \: + \: W^{-1} \: F^- \: W \: - \:
(\nabla^+ \: - \: W^{-1} \nabla^- \: W)^2 ) \right), \notag \\
G_{\alpha \beta \gamma} \: = \: 
&\frac{1}{2} \: \frac{1}{(2 \pi i)^2} \:
\int_{D^2} \sigma_{\alpha \beta \gamma}^* \: \Tr \left(
W^{-1} \: \delta W \: \wedge \: 
W^{-1} \: \delta W \: \wedge \: 
(\nabla^+ \: - \: W^{-1} \: \nabla^- W)
\right), \notag \\
G_{\alpha \beta \gamma \delta} \: = \: 
& - \: \frac{1}{6} \: \frac{1}{(2 \pi i)^2} \:
\int_{D^3} \sigma_{\alpha \beta \gamma \delta}^* \: \Tr \left(
W^{-1} \: \delta W \: \wedge \: 
W^{-1} \: \delta W \: \wedge \: 
W^{-1} \: \delta W \right),  \notag \\
G_{\alpha \beta \gamma \delta \epsilon} \: = \: 
& - \: \frac{1}{6} \: \frac{1}{(2 \pi i)^2} \:
\int_{S^3} \mu_{\alpha \beta \gamma \delta \epsilon}^* \: \Tr \left(
W^{-1} \: \delta W \: \wedge \: 
W^{-1} \: \delta W \: \wedge \: 
W^{-1} \: \delta W \right).  \notag
\end{align}
Then $G_\alpha \in \Omega^3(U_\alpha)$,
$G_{\alpha \beta} \in \Omega^2(U_\alpha \cap U_\beta)$,
$G_{\alpha \beta \gamma} \in \Omega^1(U_\alpha \cap U_\beta \cap U_\gamma)$
and
$G_{\alpha \beta \gamma \delta} \in \Omega^0(U_\alpha \cap U_\beta \cap
U_\gamma \cap U_\delta)$. 
Also, $G_{\alpha \beta \gamma \delta \epsilon}$ is an integer-valued 
continuous function
on $U_\alpha \cap U_\beta \cap
U_\gamma \cap U_\delta \cap U_\epsilon$.
The meaning of $G_{\alpha \beta}$, for example, is
\begin{align} \label{eqn5.32}
G_{\alpha \beta} \: = \:  \frac{1}{2} \: \frac{1}{(2 \pi i)^2} \:
\int_0^1 dt \: \Tr \left( \right.
& W(t)^{-1} \:  \frac{dW}{dt} \\
& \left. (F^+ \: + \: W(t)^{-1} \: F^- \: W(t) \: - \:
(\nabla^+ \: - \: W(t)^{-1} \nabla^- \: W(t))^2 ) \right). \notag
\end{align}

From (\ref{eqn5.26}),
\begin{align} \label{eqn5.33}
d G_{\alpha \beta} \: = \: & G_\beta \big|_{U_\alpha \cap U_\beta} \: - \: 
G_\alpha \big|_{U_\alpha \cap U_\beta}, \\
d G_{\alpha \beta \gamma} \: = \: & G_{\beta \gamma}
\big|_{U_\alpha \cap U_\beta \cap U_\gamma} \: - \:
G_{\alpha \gamma}
\big|_{U_\alpha \cap U_\beta \cap U_\gamma} \: + \: 
G_{\alpha \beta}
\big|_{U_\alpha \cap U_\beta \cap U_\gamma}, \notag \\
d G_{\alpha \beta \gamma \delta} \: = \: & G_{\beta \gamma \delta}
\big|_{U_\alpha \cap U_\beta \cap U_\gamma \cap U_\delta} \: - \: 
G_{\alpha \gamma \delta}
\big|_{U_\alpha \cap U_\beta \cap U_\gamma \cap U_\delta} \: + \: 
G_{\alpha \beta \delta}
\big|_{U_\alpha \cap U_\beta \cap U_\gamma \cap U_\delta} \: - \: 
G_{\alpha \beta \gamma}
\big|_{U_\alpha \cap U_\beta \cap U_\gamma \cap U_\delta}. \notag
\end{align}
Also,
\begin{align} \label{eqn5.34}
G_{\alpha \beta \gamma \delta \epsilon} \: = \: & G_{\beta \gamma \delta
\epsilon}
\big|_{U_\alpha \cap U_\beta \cap U_\gamma \cap U_\delta
\cap U_\epsilon} \: - \:
G_{\alpha \gamma \delta
\epsilon}
\big|_{U_\alpha \cap U_\beta \cap U_\gamma \cap U_\delta
\cap U_\epsilon} \: + \:G_{\alpha \beta \delta
\epsilon}
\big|_{U_\alpha \cap U_\beta \cap U_\gamma \cap U_\delta
\cap U_\epsilon} \: -  \notag \\
& G_{\alpha \beta \gamma 
\epsilon}
\big|_{U_\alpha \cap U_\beta \cap U_\gamma \cap U_\delta
\cap U_\epsilon} \: + \:G_{\alpha \beta \gamma \delta}
\big|_{U_\alpha \cap U_\beta \cap U_\gamma \cap U_\delta
\cap U_\epsilon}.
\end{align}
Thus 
\begin{equation} \label{eqn5.35}
{\cal C} \: = \: \left( G_\alpha, G_{\alpha \beta}, 
G_{\alpha \beta \gamma}, G_{\alpha \beta \gamma \delta},
G_{\alpha \beta \gamma \delta \epsilon} \right)
\end{equation}
is a $4$-cocycle for the \v{C}ech-cohomology of the complex of sheaves
\begin{equation} \label{eqn5.36}
\Z \longrightarrow \Omega^0 \longrightarrow 
\Omega^1 \longrightarrow 
\Omega^2 \longrightarrow 
\Omega^3
\end{equation}
on $B$.

Now suppose that $\left( \sigma_\alpha^\prime, \sigma_{\alpha \beta}^\prime, 
\sigma_{\alpha \beta \gamma}^\prime, 
\sigma_{\alpha \beta \gamma \delta}^\prime \right)$ is another choice of maps.
Let ${\cal C}^\prime$ be the ensuing cocycle as in (\ref{eqn5.35}).
We assume that there is a smooth $1$-parameter family of maps
$\{\sigma_\alpha(t), \sigma_{\alpha \beta}(t), 
\sigma_{\alpha \beta \gamma}(t)\}_{t \in [0,1]}$ so that for each $t
\in [0,1]$, 
$\left( \sigma_\alpha(t), \sigma_{\alpha \beta}(t), 
\sigma_{\alpha \beta \gamma}(t) \right)$ satisfies (\ref{eqn5.27}) and 
(\ref{eqn5.28}),
\begin{equation} \label{eqn5.37}
\left( \sigma_\alpha(0), \sigma_{\alpha \beta}(0), 
\sigma_{\alpha \beta \gamma}(0) \right) \: = \:
\left( \sigma_\alpha, \sigma_{\alpha \beta}, 
\sigma_{\alpha \beta \gamma} \right)
\end{equation}
and
\begin{equation} \label{eqn5.38}
\left( \sigma_\alpha(1), \sigma_{\alpha \beta}(1), 
\sigma_{\alpha \beta \gamma}(1) \right) \: = \:
\left( \sigma^\prime_\alpha, \sigma^\prime_{\alpha \beta}, 
\sigma^\prime_{\alpha \beta \gamma} \right).
\end{equation}
We do not assume that the homotopies from $\sigma_{\alpha \beta \gamma}$ to
$\sigma^\prime_{\alpha \beta \gamma}$ extend to a homotopy from
$\sigma_{\alpha \beta \gamma \delta}$ to
$\sigma^\prime_{\alpha \beta \gamma \delta}$.

Define $\Sigma_\alpha \: : \: [0, 1] \times D^0 \rightarrow {\cal I}_\alpha$
by $\Sigma_\alpha(t,x) \: = \: (\sigma_\alpha(t))(x)$, and similarly for
$\Sigma_{\alpha \beta}
 \: : \: [0, 1] \times D^1 \rightarrow {\cal I}_{\alpha \beta}$ and 
$\Sigma_{\alpha \beta \gamma}
 \: : \: [0, 1] \times D^2 \rightarrow {\cal I}_{\alpha \beta \gamma}$. 
Put $\Psi_{\alpha \beta \gamma \delta} \: = \:
- \: \Sigma_{\beta \gamma \delta} \: + \: 
\Sigma_{\alpha \gamma \delta} \: - \: 
\Sigma_{\alpha \beta \delta} \: + \: 
\Sigma_{\alpha \beta \gamma} \: + \: 
\sigma^\prime_{\alpha \beta \gamma \delta} \: - \:
\sigma_{\alpha \beta \gamma \delta}$, a map from
$S^3$ to 
${\cal I}_{\alpha \beta \gamma \delta}$.

Put
\begin{align} \label{eqn5.39}
H_{\alpha} \: = \: & \frac{1}{2} \: \frac{1}{(2 \pi i)^2} \:
\int_{[0,1] \times D^0} \Sigma_{\alpha}^* \: \Tr \left(
W^{-1} \: \delta W \: \wedge \: (F^+ \: + \: W^{-1} \: F^- \: W \: - \:
(\nabla^+ \: - \: W^{-1} \nabla^- \: W)^2 ) \right), \notag \\
H_{\alpha \beta} \: = \: 
& \frac{1}{2} \: \frac{1}{(2 \pi i)^2} \:
\int_{[0,1] \times D^1} \Sigma_{\alpha \beta}^* \: \Tr \left(
W^{-1} \: \delta W \: \wedge \: 
W^{-1} \: \delta W \: \wedge \: 
(\nabla^+ \: - \: W^{-1} \: \nabla^)
\right), \notag \\
H_{\alpha \beta \gamma} \: = \: 
& - \: \frac{1}{6} \: \frac{1}{(2 \pi i)^2} \:
\int_{[0,1] \times D^2} \Sigma_{\alpha \beta \gamma}^* \: \Tr \left(
W^{-1} \: \delta W \: \wedge \: 
W^{-1} \: \delta W \: \wedge \: 
W^{-1} \: \delta W \right),  \notag \\
H_{\alpha \beta \gamma \delta} \: = \: 
& - \: \frac{1}{6} \: \frac{1}{(2 \pi i)^2} \:
\int_{S^3} \Psi_{\alpha \beta \gamma \delta}^* \: \Tr \left(
W^{-1} \: \delta W \: \wedge \: 
W^{-1} \: \delta W \: \wedge \: 
W^{-1} \: \delta W \right).  \notag
\end{align}
Then $H_\alpha \in \Omega^2(U_\alpha)$, 
$H_{\alpha \beta} \in \Omega^1(U_\alpha \cap U_\beta)$, 
$H_{\alpha \beta \gamma} \in \Omega^0(U_\alpha \cap U_\beta
\cap U_\gamma)$,
$H_{\alpha \beta \gamma \delta}$ is an integer-valued continuous function
on ${\cal I}_{\alpha \beta \gamma \delta}$ and
\begin{align} \label{eqn5.40}
G_\alpha^\prime \: - \: G_\alpha \: & = \: d H_\alpha, \\
G_{\alpha \beta}^\prime \: - \: G_{\alpha \beta} \: & = \: d 
H_{\alpha \beta} \: + \: H_\beta \: - \: H_\alpha, \notag \\
G_{\alpha \beta \gamma}^\prime \: - \: G_{\alpha \beta \gamma} \: & = \: d 
H_{\alpha \beta \gamma} \: + \: H_{\beta \gamma} \: - \: 
H_{\alpha \gamma} \: + \: H_{\alpha \beta}, \notag \\
G_{\alpha \beta \gamma \delta}^\prime \: - \: 
G_{\alpha \beta \gamma \delta} \: & = \:  
H_{\alpha \beta \gamma \delta} \: + \: H_{\beta \gamma \delta} \: - \: 
H_{\alpha \gamma \delta} \: + \: H_{\alpha \beta \delta} \: - \:
H_{\alpha \beta \gamma}, \notag \\
G_{\alpha \beta \gamma \delta \epsilon}^\prime \: - \: 
G_{\alpha \beta \gamma \delta \epsilon} \: & = \:  
H_{\beta \gamma \delta \epsilon} \: - \: H_{\alpha \gamma \delta
\epsilon} \: + \: 
H_{\alpha \beta \delta \epsilon} \: - \: 
H_{\alpha \beta \gamma \epsilon} \: + \:
H_{\alpha \beta \gamma \delta}. \notag
\end{align}
In other words, ${\cal C}^\prime \: - \: {\cal C}$ is the coboundary of
the $3$-cochain
\begin{equation} \label{eqn5.41}
\left( H_\alpha, H_{\alpha \beta}, H_{\alpha \beta \gamma},
H_{\alpha \beta \gamma \delta}  
\right).
\end{equation}

In summary, our input data consisted of the maps
$\left( \sigma_\alpha, \sigma_{\alpha \beta}, \sigma_{\alpha \beta \gamma}
\right)$ satisfying (\ref{eqn5.27}) and (\ref{eqn5.28}),
defined up to homotopy, with the property that if 
$U_\alpha \cap U_\beta \cap U_\gamma \cap U_\delta
\neq \emptyset$ then the map from $S^2$ to 
${\cal I}_{\alpha \beta \gamma \delta}$, given by
$\sigma_{\beta \gamma \delta} 
\big|_{U_\alpha \cap U_\beta \cap U_\gamma \cap U_\delta} \: - \:
\sigma_{\alpha \gamma \delta} 
\big|_{U_\alpha \cap U_\beta \cap U_\gamma \cap U_\delta} \: + \:
\sigma_{\alpha \beta \delta} 
\big|_{U_\alpha \cap U_\beta \cap U_\gamma \cap U_\delta} \: - \:
\sigma_{\alpha \beta \gamma} 
\big|_{U_\alpha \cap U_\beta \cap U_\gamma \cap U_\delta}$,
extends to a map from $D^3$ to ${\cal I}_{\alpha \beta \gamma \delta}$.
From this we obtained a 
Deligne cohomology class on $B$
of degree $4$ \cite[Chapter 1.5]{Brylinski (1993)}.
From (\ref{eqn5.9}), its ``curvature'' is the $4$-form
\begin{equation} \label{eqn5.42}
\left( \int_Z \widehat{A} \left(R^{TZ}/2\pi i \right)
\: \wedge \: \ch \left( F^V/2\pi i \right)
 \right)^{(4)} \in
\Omega^4(B).
\end{equation}

\subsubsection{$\HH^{2k}$} \label{HH2k}

Let us write the degree-$(2k \: - \: 1)$ component $\CS^{2k-1}$ of $\CS$ as
$\CS^{2k-1} \: = \: \sum_{l=0}^{2k-1} \CS^{l, 2k-l-1}$, with
$\CS^{l, 2k-l-1} \in \Omega^l({\cal I}) \: \widehat{\otimes} \:
\Omega^{2k-l-1}(B)$. Given $\{U_\alpha\}_{\alpha \in I}$, we suppose that
for $0 \: \le \: l \: \le \: 2k \: - \: 1$ there exist maps 
$\sigma_{\alpha_0 \ldots \alpha_l} \: : \: D^l \rightarrow
{\cal I}_{\alpha_0 \ldots \alpha_l}$ such that
\begin{equation} \label{eqn5.43}
\partial \sigma_{\alpha_0 \ldots \alpha_l} \: = \:
\sum_{m=0}^l \: (-1)^m \: \sigma_{\alpha_0 \ldots
\widehat{\alpha_m} \ldots \alpha_l}.
\end{equation}
If $U_{\alpha_0} \cap \ldots \cap U_{\alpha_{2k}} \neq \emptyset$, define
$\mu_{{\alpha_0} \ldots {\alpha_{2k}}} \: : \: S^{2k-1} \rightarrow
{\cal I}_{{\alpha_0} \ldots {\alpha_{2k}}}$ by
\begin{equation} \label{eqn5.44}
\mu_{{\alpha_0} \ldots {\alpha_{2k}}} \: = \:
\sum_{m=0}^{2k} \: (-1)^m \: 
\sigma_{\alpha_0 \ldots
\widehat{\alpha_m} \ldots \alpha_{2k}} 
\big|_{U_{\alpha_0} \cap \ldots \cap U_{\alpha_{2k}}}.
\end{equation}
Again, we assume that the map $\mu_{{\alpha_0} \ldots {\alpha_{2k}}}$ has
been parametrized so as to be smooth.  
Put 
$G_{\alpha_0} \: = \: \widehat{\eta}^{(2k-1)} \: + \:
\CS^{0, 2k-1}_{\alpha_0}$. 
For $1 \: \le \: l \le \: 2k \: - \: 1$,
put
\begin{equation} \label{eqn5.45}
G_{\alpha_0 \ldots \alpha_l} \: = \:
\int_{D^l} \sigma_{\alpha_0 \ldots \alpha_l}^* \: 
\CS^{l, 2k-l-1},
\end{equation}
and put
\begin{equation} \label{eqn5.46}
G_{\alpha_0 \ldots \alpha_{2k}} \: = \:
\int_{S^{2k-1}} \mu_{\alpha_0 \ldots \alpha_{2k}}^* \: 
\CS^{2k-1, 0},
\end{equation}
Then $G_{\alpha_0 \ldots \alpha_l} \in 
\Omega^{2k-l-1}(U_{\alpha_0} \cap \ldots \cap U_{\alpha_l})$.

\begin{lemma}
$G_{\alpha_0 \ldots \alpha_{2k}}$ is an integer-valued
continous function on 
$U_{\alpha_0} \cap \ldots \cap U_{\alpha_{2k}}$.
\end{lemma}
\begin{pf}
From (\ref{eqn5.10}) and (\ref{eqn5.12}),
\begin{align}
\CS^{2k-1, 0} \: & = \:  \frac{1}{(2\pi i)^k \: (k-1)!} \:
\Tr \left( W^{-1} \: \delta W \right)^{2k-1} \: \int_0^1 t^{k-1} \:
(1-t)^{k-1} \: dt \\
& = \:  \frac{1}{(2\pi i)^k} \:
\frac{(k-1)!}{(2k-1)!} \:
\Tr \left( W^{-1} \: \delta W \right)^{2k-1}. \notag
\end{align}
Thus 
\begin{equation}
G_{\alpha_0 \ldots \alpha_{2k}} \: = \:
\frac{1}{(2\pi i)^k} \:
\frac{(k-1)!}{(2k-1)!} \: 
\int_{S^{2k-1}} \mu_{\alpha_0 \ldots \alpha_{2k}}^* \: 
\Tr \left( W^{-1} \: \delta W \right)^{2k-1}.
\end{equation}
From \cite[p. 237]{Bott-Seeley (1978)}, this is integer-valued. 
For this 
to be true, it is important that we are integrating over $S^{2k-1}$ and
not over an arbitrary $(2k \: - \: 1)$-dimensional manifold.
\end{pf}

We have
\begin{equation} \label{eqn5.47}
dG_{\alpha_0} \: = \: \left( \int_Z \widehat{A} \left(R^{TZ}/2\pi i \right)
\: \wedge \: \ch \left( F^V/2\pi i \right)
 \right)^{(2k)} \in
\Omega^{2k}(B)
\end{equation}
and for $1 \: \le \: l \le \: 2k \: - \: 1$
\begin{equation} \label{eqn5.48}
dG_{\alpha_0 \ldots \alpha_l} \: = \: 
\sum_{m=0}^l \: (-1)^m \: G_{\alpha_0 \ldots \widehat{\alpha_m} \ldots 
\alpha_l} \big|_{U_{\alpha_0} \cap \ldots \cap U_{\alpha_l}}.
\end{equation}
Also,

\begin{equation} \label{eqn5.49}
G_{\alpha_0 \ldots \alpha_{2k}} \: = \: 
\sum_{m=0}^{2k} \: (-1)^m \: G_{\alpha_0 \ldots \widehat{\alpha_m} \ldots 
\alpha_{2k}} \big|_{U_{\alpha_0} \cap \ldots \cap U_{\alpha_{2k}}}.
\end{equation}
Thus
\begin{equation} \label{eqn5.50}
{\cal C} \: = \: \left( G_{\alpha_0}, \ldots, 
G_{\alpha_0 \ldots \alpha_{2k}} \right)
\end{equation}
is a $2k$-cocycle for the \v{C}ech-cohomology of the complex of sheaves
\begin{equation} \label{eqn5.51}
\Z \longrightarrow \Omega^0 \longrightarrow \ldots \longrightarrow
\Omega^{2k-1}
\end{equation}
on $B$.

Now suppose that $\left( \sigma_{\alpha_0}^\prime, \ldots, 
\sigma_{\alpha_0 \ldots \alpha_{2k-1}}^\prime \right)$
is another choice of maps.
Let ${\cal C}^\prime$ be the ensuing cocycle as in (\ref{eqn5.50}).
We assume that for $0 \: \le \: l \: \le \: 2k-2$
there is a smooth $1$-parameter family of maps
$\{ \sigma_{\alpha_0 \ldots \alpha_{l}}(t)\}_{t \in [0,1]}$ 
so that for each $t \in [0,1]$, 
$\left( \sigma_{\alpha_0}(t), \ldots, 
\sigma_{\alpha_0 \ldots \alpha_{l}}(t) \right)$
 satisfies (\ref{eqn5.43}),
\begin{equation} \label{eqn5.52}
\sigma_{\alpha_0 \ldots \alpha_{l}}(0) \: = \: 
\sigma_{\alpha_0 \ldots \alpha_{l}}
\end{equation}
and
\begin{equation} \label{eqn5.53}
\sigma_{\alpha_0 \ldots \alpha_{l}}(1) \: = \: 
\sigma_{\alpha_0 \ldots \alpha_{l}}^\prime.
\end{equation}
We do not assume that the homotopies from $\sigma_{\alpha_0 \ldots
\alpha_{2k-2}}$ to
$\sigma^\prime_{\alpha_0 \ldots
\alpha_{2k-2}}$ extend to a homotopy from
$\sigma_{\alpha_0 \ldots \alpha_{2k-1}}$ to
$\sigma^\prime_{\alpha_0 \ldots \alpha_{2k-1}}$.

\begin{proposition}
Under these assumptions, ${\cal C}$ and ${\cal C}^\prime$ are
cohomologous.
\end{proposition}
\begin{pf}
For $0 \: \le \: l \: \le \: 2k \: - \: 2$,
define 
$\Sigma_{\alpha_0 \ldots \alpha_l} \: : \: [0, 1] \times D^l 
\rightarrow {\cal I}_{\alpha_0 \ldots \alpha_l}$
by $\Sigma_{\alpha_0 \ldots \alpha_l}(t,x) \: = \: 
(\sigma_{\alpha_0 \ldots \alpha_l}(t))(x)$.
Put 
\begin{equation} \label{eqn5.54}
\Psi_{\alpha_0 \ldots \alpha_{2k-1}} \: = \:
\sigma^\prime_{\alpha_0 \ldots \alpha_{2k-1}} \: - \:
\sigma_{\alpha_0 \ldots \alpha_{2k-1}} \: - \:
\sum_{m=0}^{2k-1} (-1)^m \: 
\Sigma_{\alpha_0 \ldots \widehat{\alpha_m} \ldots \alpha_{2k-1}},
\end{equation}
a map from
$S^{2k-1}$ to 
${\cal I}_{\alpha_0 \ldots \alpha_{2k-1}}$.

For $0 \: \le \: l \: \le 2k \: - \: 2$, put 
\begin{equation} \label{eqn5.55}
H_{\alpha_0 \ldots \alpha_l} \: = \: 
\int_{[0,1] \times D^l} \Sigma_{\alpha_0 \ldots \alpha_l}^* \: 
\CS^{l+1, 2k - l - 2}.
\end{equation}
Then $H_{\alpha_0 \ldots \alpha_l} \in 
\Omega^{2k-l-2}(U_{\alpha_0} \cap \ldots \cap
U_{\alpha_l})$.
Put
\begin{equation} \label{eqn5.56}
H_{\alpha_0 \ldots \alpha_{2k-1}} \: = \: 
\int_{S^{2k-1}} \Psi_{\alpha_0 \ldots \alpha_{2k-1}}^* \: 
\CS^{2k-1,0}.
\end{equation}
Then
$H_{\alpha_0 \ldots \alpha_{2k-1}}$ is an integer-valued continuous function
on ${\cal I}_{\alpha_0 \ldots \alpha_{2k-1}}$. Furthermore,
for $0 \: \le \: l \: \le \: 2k \: - \: 2$, 
\begin{equation} \label{eqn5.57}
G_{\alpha_0 \ldots \alpha_l}^\prime \: - \: 
G_{\alpha_0 \ldots \alpha_l} \: = \: d 
H_{\alpha_0 \ldots \alpha_l} \: + \: 
\sum_{m=0}^l \: (-1)^m \:
H_{\alpha_0 \ldots \widehat{\alpha_m} \ldots \alpha_l},
\end{equation}
\begin{equation} \label{eqn5.58}
G_{\alpha_0 \ldots \alpha_{2k-1}}^\prime \: - \: 
G_{\alpha_0 \ldots \alpha_{2k-1}} \: = \:
H_{\alpha_0 \ldots \alpha_{2k-1}} \: + \: 
\sum_{m=0}^{2k-1} \: (-1)^m \:
H_{\alpha_0 \ldots \widehat{\alpha_m} \ldots \alpha_{2k-1}}
\end{equation}
and
\begin{equation} \label{eqn5.59}
G_{\alpha_0 \ldots \alpha_{2k}}^\prime \: - \: 
G_{\alpha_0 \ldots \alpha_{2k}} \: = \:
\sum_{m=0}^{2k} \: (-1)^m \:
H_{\alpha_0 \ldots \widehat{\alpha_m} \ldots \alpha_{2k}}.
\end{equation}
In other words, ${\cal C}^\prime \: - \: {\cal C}$ is the coboundary of
the $(2k \: - \: 1)$-cochain
\begin{equation} \label{eqn5.60}
\left( H_{\alpha_0}, \ldots, 
H_{\alpha_0 \ldots \alpha_{2k-1}}  
\right).
\end{equation}
Thus ${\cal C}$ and ${\cal C}^\prime$ are cohomologous.
\end{pf}

In summary, our input data consisted of the maps
$\left( \sigma_{\alpha_0}, \ldots,  \sigma_{\alpha_0 \ldots \alpha_{2k-2}}
\right)$ satisfying (\ref{eqn5.43}),
defined up to homotopy, with the property that if 
$U_{\alpha_0} \cap \ldots \cap U_{\alpha_{2k-1}}
\neq \emptyset$ then the map from $S^{2k-2}$ to 
${\cal I}_{\alpha_0 \ldots \alpha_{2k-1}}$, given by
$\sum_{m=0}^{2k-1} \: (-1)^m \: \sigma_{\alpha_0 \ldots \widehat{\alpha_m}
\ldots \alpha_{2k-1}} 
\big|_{U_{\alpha_0} \cap \ldots \cap U_{\alpha_{2k-1}}}$,
extends to a map from $D^{2k-1}$ to 
${\cal I}_{\alpha_0 \ldots \alpha_{2k-1}}$.
From this we obtained a 
Deligne cohomology class on $B$
of degree $2k$ \cite[Chapter 1.5]{Brylinski (1993)}.
From (\ref{eqn5.9}), its ``curvature'' is the $2k$-form
\begin{equation} \label{eqn5.61}
\left( \int_Z \widehat{A} \left(R^{TZ}/2\pi i \right)
\: \wedge \: \ch \left( F^V/2\pi i \right)
 \right)^{(2k)} \in
\Omega^{2k}(B).
\end{equation}

\subsubsection{Topological interpretation}

Let us note that we can always add a trivial vector bundle
$B \times \C^N$, with a trivial connection, to both 
$\Ker(D)_+$ and $\Ker(D)_-$ and carry out the preceding constructions for this
stabilized vector bundle.
Thus it is only the stabilized class of
$\Ker(D)$ that matters.

Recall that the nerve ${\cal N}$ of a covering $\{U_\alpha\}_{\alpha \in I}$
is a certain simplicial complex which has one $k$-simplex for each
nonempty intersection $U_{\alpha_0} \cap \ldots \cap U_{\alpha_k}$.
Let ${\cal N}^{(k)}$ denote the $k$-skeleton of ${\cal N}$.

Consider the space $X$ obtained by gluing together
$\{(U_{\alpha_0} \cap \ldots \cap U_{\alpha_l}) \times \Delta^l 
\}_{l=0}^\infty$ using the embeddings
$(U_{\alpha_0} \cap \ldots \cap U_{\alpha_l}) \times \partial_m \Delta^l
\rightarrow
(U_{\alpha_0} \cap \ldots \cap \widehat{U_{\alpha_m}} \cap \ldots
 \cap U_{\alpha_l}) \times \Delta^l$. 
There is a continuous map $\rho \: : \: X \rightarrow
{\cal N}$ which contracts each 
$U_{\alpha_0} \cap \ldots \cap U_{\alpha_l}$ to a point.
We now assume that $\{U_\alpha\}_{\alpha \in I}$ is a good covering,
meaning that each 
$U_{\alpha_0} \cap \ldots \cap U_{\alpha_l}$ is contractible.
Then each preimage of $\rho$ is contractible and in our case it follows that
$\rho$ is a homotopy equivalence
\cite{Dold-Thom (1958)}. 

There is an obvious $\Z_2$-graded vector bundle $V$ on $X$ whose restriction to
$(U_{\alpha_0} \cap \ldots \cap U_{\alpha_l}) \times \Delta^l$ pulls back
from $\Ker(D) \big|_{U_{\alpha_0} \cap \ldots \cap U_{\alpha_l}}$.
Suppose that we have the maps $\{\sigma_{\alpha_0 \ldots \alpha_l}
\}_{l=0}^{2k-1}$ of the previous subsubsection.  Then these isometries show
that $V$ is trivial as a $\Z_2$-graded vector bundle 
on $\rho^{-1} \left( {\cal N}^{(2k-1)} \right)$, i.e. that there is an
isomorphism there from $V_+$ to $V_-$.

Recall that there is a filtration
$\KK^*(X) = 
\KK^*_0(X) \supset \KK^*_1(X) \supset \ldots$ of $\KK$-theory, where
$\KK^*_i(X)$ consists of the elements 
$x$ of ${\KK}^*(X)$
with the property that for any finite simplicial complex $Y$ 
of dimension less than
$i$ and any continuous map $f \: : \: Y \rightarrow X$, 
$f^* x \: = \: 0$
\cite[Section 2]{Atiyah-Hirzebruch (1961)}. By definition, the 
filtration is homotopy-invariant.
It gives rise
to the Atiyah-Hirzebruch spectral sequence to compute $\KK^*(X)$, with
$E_2$-term $E_2^{p,q} \: = \: \HH^p(X; \KK^q(\pt.))$ and
$E_\infty$-term $E_\infty^{p,q} \: = \: \KK^{p+q}_p(X)/\KK^{p+q}_{p+1}(X)$.

Using the homotopy equivalence given by $\rho$ from $X$ to the simplicial
complex ${\cal N}$, 
we see that the $\KK$-theory class of the $\Z_2$-graded vector
bundle $V$ on $X$ lies in $\KK^{0}_{2k}(X)$. Equivalently, the 
$\KK$-theory class of the $\Z_2$-graded vector
bundle $\Ker(D)$ on $B$ lies in $\KK^{0}_{2k}(B)$

Conversely, suppose that the 
$\KK$-theory class of the $\Z_2$-graded vector
bundle $\Ker(D)$ on $B$ lies in $\KK^{0}_{2k}(B)$.
Then $V$ has its $\KK$-theory class in 
$\KK^{0}_{2k}(X)$.  This means that
$V \big|_{\rho^{-1} \left( {\cal N}^{(2k-1)} \right)}$ vanishes in
$\KK^{0} \left( \rho^{-1} \left( {\cal N}^{(2k-1)} \right) \right)$.
After possibly stabilizing by trivial bundles, so that
$V_+ \big|_{\rho^{-1} \left( {\cal N}^{(2k-1)} \right)}$ is isomorphic to
$V_- \big|_{\rho^{-1} \left( {\cal N}^{(2k-1)} \right)}$, 
we obtain the existence of
the maps $\{\sigma_{\alpha_0 \ldots \alpha_l}
\}_{l=0}^{2k-1}$ of the previous subsection.

As a further point, if $\rk(\Ker(D)_+) \neq
\rk(\Ker(D)_-)$ originally then 
after adding trivial vector bundles, we may assume that
$\rk(\Ker(D)_+) \: = \: \rk(\Ker(D)_-)$. Thus it is enough to only
consider the image of the $\KK$-theory class of $\Ker(D)$ in the 
reduced $\KK$-theory group
$\widetilde{\KK}^{0}(B)$.

The Deligne cohomology class depends, {\it a priori}, on the choice of
$\{\sigma_{\alpha_0 \ldots \alpha_l}
\}_{l=0}^{2k-2}$. We now examine how many such choices there are. We
implicitly
stabilize the vector bundle $V$.
Let $\BU \langle p-1 \rangle$ 
denote the $(p-2)$-connected space which appears in the
Whitehead tower of $\Z \times \BU$. Then $[{\cal N}, \BU\langle p-1 \rangle]$, 
the set of homotopy classes of continuous maps 
from ${\cal N}$ to $\BU\langle p-1 \rangle$, 
is isomorphic to the homotopy classes of
vector bundles which are
trivial on ${\cal N}^{(p-2)}$ and which are trivialized on 
${\cal N}^{(p-3)}$. Suppose that we have such a map $\phi \: : \: 
{\cal N} \rightarrow \BU\langle p-1 \rangle$. 
The obstruction to lifting $\phi$, 
with respect to the map
$\BU\langle p \rangle 
\rightarrow \BU\langle p-1 \rangle$, 
to a map  $ \widetilde{\phi} \: : \:
{\cal N} \rightarrow \BU\langle p \rangle$ 
corresponds to the obstruction to trivializing
the vector bundle over ${\cal N}^{(p-1)}$. That is
$\widetilde{\KK}^0_p({\cal N})$ is the same as the elements of
$[{\cal N}, \BU]$ which lift to $[{\cal N}, \BU\langle p \rangle ]$.
The number of such liftings corresponds to the number of trivializations
of the vector bundle on ${\cal N}^{(p-2)}$. Let ${\cal F}_p$ denote the
homotopy fiber of the map
$\BU\langle p \rangle \rightarrow \BU\langle p-1 \rangle$. As 
$\BU\langle p \rangle \rightarrow \BU\langle p-1 \rangle$ 
is a principal fibration, the number of
liftings  to
$[{\cal N}, \BU\langle p \rangle ]$ 
of a liftable element of $[{\cal N}, \BU\langle p-1 \rangle ]$ is given by 
\begin{equation}
[{\cal N}, {\cal F}_p] \: = \: 
 [{\cal N}, \KK(\pi_{p-1}(\BU), p-2)] \: = \:
\HH^{p-2}(B; \pi_{p-1}(\BU)).
\end{equation} 
Proceeding inductively over the skeleta, we see that the
set of possible degree-$2k$ Deligne cohomology classes that we can construct
is acted upon transitively by $\bigoplus_{j=1}^\infty \HH^{2k-1-2j}(B; \Z)$.

Thus under the assumption that $\Ker(D)$ is a vector bundle on $B$
whose reduced $\KK$-theory class lies in $\widetilde{\KK}^0_{2k}(B)$,
we obtain a set of Deligne cohomology classes of degree $2k$.
 
If $B$ happens to be $(2k \: - \: 2)$-connected then there is a unique
lifting of an element of $[B, \BU]$ to $[B, \BU \langle 2k \rangle]$,
so there is no obstruction to the existence of the Deligne
cohomology class and we obtain a single such
class.

\subsection{The general case}

We no longer assume that $\Ker(D)$ forms a vector bundle on $B$. We will
essentially
reduce to the case of vector bundle kernel by means of the method of
\cite[Section 5]{Lott (1994)}.

From a general result in index theory, there are smooth finite-dimensional
subbundles $K_\pm$ of $(\pi_* E)_\pm$ and complementary subbundles
$G_\pm$ such that $D$ is diagonal with respect to the decomposition
$(\pi_* E)_\pm \: = \: G_\pm \: \oplus \: K_\pm$ and writing
$D \: = \: D_G \: + \: D_K$, in addition
$(D_G)_\pm : C^\infty(B; G_\pm) \rightarrow C^\infty(B; G_\mp)$ is fiberwise
$L^2$-invertible
\cite{Miscenko-Fomenko (1979)}.
Give $K_\pm$ a Hermitian metric $h^{K_\pm}$ and compatible
connection $\nabla^{K_\pm}$. Let $F^{K_\pm}$ denote the curvature of
$\nabla^{K_\pm}$.

Put 
\begin{equation} \label{eqn6.1}
H_\pm \: = \: (\pi_* E)_\pm \: \oplus \: K_\mp \: = \: 
G_\pm \: \oplus \: K_\pm \: \oplus \: K_\mp,
\end{equation}
where the factor $K_\mp$ has the metric $h^{K_\mp}$ and
connection $\nabla^{K_\mp}$.
Let $\phi \: : \: [0, \infty) \rightarrow [0,1]$ be a smooth bump function
such that there exist $\delta, \Delta \: > \: 0$ satisfying
\begin{equation} \label{eqn6.2}
\phi(t) \: = \:
\begin{cases}
0 & \text{ if } t \in (0, \delta), \\
1 & \text{ if } t > \Delta.
\end{cases}
\end{equation}
For $s \in \R$, define $R_+(s) \: : \: C^\infty(B; H_+) \rightarrow 
C^\infty(B; H_-)$ by
\begin{equation} \label{eqn6.3}
R_+(s) \: = \:
\begin{pmatrix}
0 & 0 & 0 \\
0 & 0 & \alpha \phi(s) \\
0 & \alpha \phi(s) & 0
\end{pmatrix}
\end{equation}
and let $R_-(s)$ be its adjoint.
Define a family ${\cal A}_s$ of superconnections on $H$ by
\begin{equation} \label{eqn6.4}
{\cal A}_s \: = \: (A_s \: \oplus \: \nabla^K) \: + \: s \:
R(s).
\end{equation}
If $s \in (0, \delta)$ then ${\cal A}_s \: = \: A_s \: \oplus \: \nabla^K$, 
while if
$s \: > \: \Delta$ then the component of ${\cal A}_s$ of degree zero,
with respect to $B$, which maps $C^\infty(B; H_+)$ to 
$C^\infty(B; H_-)$ is $s$ times
\begin{equation} \label{eqn6.5}
\begin{pmatrix}
(D_G)_+ & 0 & 0 \\
0 & (D_K)_+ & \alpha \\
0 & \alpha & 0
\end{pmatrix}.
\end{equation}
If $\alpha$ is sufficiently large, which we will assume,
then the operator in (\ref{eqn6.5}) is
$L^2$-invertible.

Define $\widetilde{\eta}$ and $\widehat{\eta}$ as in (\ref{eqn5.3}) and 
(\ref{eqn5.5}),
using ${\cal A}_s$ instead of $A_s$; the idea of this sort of
$s$-dependent definition of ${\cal A}_s$ is taken from
\cite{Melrose-Piazza (1997)}. Then using the formula
\begin{equation} \label{eqn6.6}
\frac{d}{ds} \: \Tr_s \left( e^{- \: {\cal A}_s^2} \right) \: = \: - \:
d \: \Tr_s \left( \frac{d{\cal A}_s}{ds} \: e^{- \: {\cal A}_s^2} \right), 
\end{equation}
one finds
\begin{equation} \label{eqn6.7}
d \widehat{\eta} \: = \:
\int_Z \widehat{A} \left(R^{TZ}/2\pi i \right)
\: \wedge \: \ch \left( F^V/2\pi i \right)
 \: - \: \ch \left( F^{K}/ 2\pi i \right).
\end{equation}

Now suppose that the image of $\Ind(D)$ under the map
${\KK}^0(B) \rightarrow \widetilde{\KK}^0(B)$ 
lies in $\widetilde{\KK}^0_{2k}(B)$.
After possibly adding trivial bundles to $K$,
we can perform the constructions of Subsubsection \ref{HH2k},
replacing $\Ker(D)_\pm$ by $K_\pm$. For example,
\begin{equation} \label{eqn6.8}
G_{\alpha_0} \: = \: \widehat{\eta}^{2k-1} \: + \: \CS_{\alpha_0}^{0,2k-1}, 
\end{equation}
where $\widehat{\eta}$ is defined using ${\cal A}$ and 
$\CS_{\alpha_0}^{0,2k-1}$ is defined using $K_\pm$ and 
$\nabla^{K_\pm}$. In this way, we obtain an
explicit Deligne cocycle on $B$ of degree $2k$.

\begin{theorem} \label{thm2}
The Deligne cohomology class
is independent of the choices of $K$, $h^K$, $\nabla^K$, $\alpha$ and $\phi$.
Its curvature is the $2k$-form
\begin{equation}
\left( \int_Z \widehat{A} \left(R^{TZ}/2\pi i \right)
\: \wedge \: \ch \left( F^V/2\pi i \right) \right)^{(2k)} \in
\Omega^{2k}(B).
\end{equation}
\end{theorem}
\begin{pf}
Suppose first that we fix $K$, $h^K$ and $\nabla^K$. Let
$\{\alpha(\epsilon)\}_{\epsilon \in [0,1]}$ and
$\{\phi(\epsilon)\}_{\epsilon \in [0,1]}$ be
smooth $1$-parameter families.

From the $\Z_2$-graded analog of (\ref{eqn2.22}), as
the small-$s$ behavior of ${\cal A}_s$ is independent of $\epsilon$, it
follows that $\frac{d\widehat{\eta}}{d\epsilon}$ is 
exact. Thus $\widehat{\eta}(1) \: - \: \widehat{\eta}(0) \: = \: 
d{\cal S}$ for some ${\cal S} \in \Omega^{even}(B)$. Then the difference of the
Deligne cocycles defined using ($\alpha(1)$ and $\phi(1)$) vs.
($\alpha(0)$ and $\phi(0)$) is the coboundary of the cochain
$({\cal S}^{(2k-2)} \big|_{U_{\alpha_0}}, 0, \ldots, 0)$.
As any two choices of $\alpha$ and
$\phi$ can be joined by such paths, it follows that the Deligne
cohomology class is independent of the choices of $\alpha$ and $\phi$.

Now suppose that we have smooth $1$-parameter families
$\{h^K(\epsilon)\}_{\epsilon \in [0,1]}$ and
$\{\nabla^K(\epsilon)\}_{\epsilon \in [0,1]}$.
From the $\Z_2$-graded analog of (\ref{eqn2.22}), on $B$ we have
\begin{equation} \label{eqn6.9}
\frac{d\widehat{\eta}(\epsilon)}{d\epsilon} \: = \:
\frac{1}{2\pi i} \: \Tr_s \left( \frac{d\nabla^K}{d\epsilon} \:
e^{- F^K/2\pi i} \right) \pmod{\Image(d)}.
\end{equation}
On the other hand, on ${\cal I} \times B$,
\begin{equation} \label{eqn6.10}
\frac{d\CS(\epsilon)}{d\epsilon} \: = \:
- \: \frac{1}{2\pi i} \: \Tr_s \left( \frac{d\nabla^K}{d\epsilon} \:
e^{- F^K/2\pi i} \right) \pmod{\Image(\delta \: + \: d)}.
\end{equation}
It follows that on ${\cal I} \times B$,
\begin{equation} \label{eqn6.11}
\left( \widehat{\eta}^{(2k-1)} \: + \: \CS^{2k-1} \right)(1) \: - \:
\left( \widehat{\eta}^{(2k-1)} \: + \: \CS^{2k-1} \right)(0) \: = \: 
(\delta \: + \: d)  {\cal S}
\end{equation}
for some ${\cal S} \in \Omega^{even}({\cal I} \times B)$. Then with respect to 
the cocycle ${\cal C}$ of (\ref{eqn5.50}), ${\cal C}(1) \: - \: {\cal C}(0)$ is
the coboundary of the cochain
\begin{equation} \label{eqn6.12}
\left( \int_{D^0} \sigma_{\alpha_0}^* {\cal S}^{0,2k-2 }, 
\int_{D^1} \sigma_{\alpha_0, \alpha_1}^* {\cal S}^{1,2k-3 }, \ldots, 
\int_{D^{2k-2}} \sigma_{\alpha_0, \ldots, \alpha_{2k-2}}^* 
{\cal S}^{2k-2,0 }, 0
\right). 
\end{equation}

Finally, suppose that $K^\prime$ is another choice of $K$.
As both $[K_+ \: - \: K_-]$ and $[K^\prime_+ \: - \: K^\prime_-]$
represent $\Ind(D)$ in $\KK^0(B)$, there are vector
bundles $L$ and $L^\prime$ such that there are topological isomorphisms
\begin{align} \label{eqn6.13}
t_+ \: : \: K_+ \: \oplus \: L & \rightarrow 
K^\prime_+ \: \oplus \: L^\prime, \\
t_- \: : \: K_- \: \oplus \: L & \rightarrow \:  
K^\prime_- \: \oplus \: L^\prime. \notag
\end{align}
Choose Hermitian metrics $h^{L}$ and $h^{L^\prime}$, and compatible
connections $\nabla^L$ and $\nabla^{L^\prime}$. 
Put $\widetilde{K}_\pm \: = \: K_\pm \: \oplus \: L$ and
$\widetilde{K}^\prime_\pm \: = \: K^\prime_\pm \: \oplus \: L^\prime$.
Now $\widehat{\eta}$ 
is unchanged if we define it in the obvious way on
\begin{equation} \label{eqn6.14}
\widetilde{H}_\pm \: = \: 
G_\pm \: \oplus \: \widetilde{K}_\pm \: \oplus \: \widetilde{K}_\mp,
\end{equation}
instead of $H_\pm$. 
Similarly, 
$\widehat{\eta}^\prime$, the eta-form corresponding to $K^\prime$, can
be computed on
\begin{equation} \label{eqn6.15}
\widetilde{H}^\prime_\pm \: = \: 
G_\pm \: \oplus \: \widetilde{K}^\prime_\pm \: \oplus \: 
\widetilde{K}^\prime_\mp
\end{equation}
The isomorphism $t$ induces an isomorphism $\widehat{t} \: : \: 
\widetilde{H}_\pm \rightarrow \widetilde{H}^\prime_\pm$. Hence 
we can also consider 
$\widehat{\eta}^\prime$ to be computed on $\widetilde{H}_\pm$ using the
pullback superconnection $\widehat{t}^* \widetilde{\cal A}^\prime$. 
Now let us compare $\widehat{t}^* \widetilde{\cal A}^\prime$
to $\widetilde{\cal A}$.
The differences in the components of degree zero and two,
with respect to $B$, are
finite-rank and, as in the preceding arguments, the ensuing Deligne cohomology
class is unchanged.  Hence we may assume that
$\widehat{t}^* \widetilde{\cal A}^\prime$ and
$\widetilde{\cal A}$ have the same components of degree zero and two.
The difference in the degree-one components comes
from the difference between $t^* \nabla^{\widetilde{K}^\prime}$ and
$\nabla^{\widetilde{K}}$. We can apply the preceding argument 
concerning the independence with respect to the choice of connection,
to conclude that the
Deligne cocycle computed with $K$ is cohomologous to that computed with
$K^\prime$.

The curvature statement follows from (\ref{eqn6.7}).
\end{pf}

Finally, suppose that $\dim(Z)$ is odd. Consider the fiber bundle
$(S^1 \times S^1 \times M) \rightarrow (S^1 \times B)$. Give the fiber
circle a length of $1$. As in
\cite[Pf. of Theorem 2.10]{Bismut-Freed (1986)} there is a canonical
family $\widetilde{D}$
of Dirac-type operators on the new fiber bundle, whose index can be 
trivialized on $\{1\} \times B \subset S^1 \times B$.  Suppose that 
$\Ind(D)$ 
lies in ${\KK}^{1}_{2k-1}(B)$. Then the image of
$\Ind(\widetilde{D})$ under the map $\KK^{0}(S^1 \times B) \rightarrow 
\widetilde{\KK}^{0}(S^1 \times B)$
lies in $\widetilde{\KK}^{0}_{2k}(S^1 \times B)$,
and so we can construct
the corresponding Deligne cohomology classes on $S^1 \times B$ of degree $2k$.
Integrating over the circle in Deligne cohomology
\cite[Section 6.5]{Brylinski (1993)}, we obtain Deligne
cohomology classes on $B$ of degree $2k-1$. The set of Deligne cohomology
classes so obtained
is acted upon transitively by
$\bigoplus_{j=1}^\infty \HH^{2k-1-2j}(S^1 \times B, \{1 \} \times B; \Z) \:
= \: \bigoplus_{j=1}^\infty \HH^{(2k-1)-1-2j}(B; \Z)$.


\begin{thebibliography}{9}

\bibitem{Atiyah-Hirzebruch (1961)} 
M. Atiyah and F. Hirzebruch,
``Vector Bundles and Homogeneous Spaces'',
Proc. Sympos. Pure Math. Vol. III, 
American Mathematical Society, Providence,
p. 7--38 (1961)

\bibitem{Atiyah-Patodi-Singer (1975)} M. Atiyah, V. Patodi and 
I. Singer, ``Spectral Asymmetry and Riemannian Geometry I'', 
Math. Proc. Cambridge Philos. Soc. 77, p. 43-69 (1975)

\bibitem{Atiyah-Singer (1984)} M. Atiyah and I. M.  Singer,
``Dirac operators Coupled to Vector Potentials'', Proc. Nat. Acad. Sci.
U.S.A. 81, p. 2597-2600 (1984)

\bibitem{Berline-Getzler-Vergne (1992)} N. Berline, E. Getzler and M. Vergne,
\underline{Heat Kernels and the
Dirac Operator}, Grundl. der Math. Wiss. 298, Springer, 
Berlin-Heidelberg-New York (1992)

\bibitem{Bismut (1986)} J.-M. Bismut, ``The Index Theorem for Families of
Dirac Operators : Two Heat Equation Proofs'', Inv. Math. 83, p. 91-151 (1986)

\bibitem{Bismut-Cheeger (1989)} J.-M. Bismut and J. Cheeger,
``$\eta$-Invariants and Their Adiabatic Limits'', J. of the Amer. Math. Soc.
2, p. 33-70 (1989)

\bibitem{Bismut-Freed (1986)} J.-M. Bismut and D. Freed,
``The Analysis of Elliptic Families I and II'', 
Comm. Math. Phys. 106, p. 159-176 and 
107, p. 103-163 (1986)

\bibitem{Bott-Seeley (1978)} R. Bott and R. Seeley, 
``Some remarks on the paper of Callias'',
Comm. Math. Phys. 62, p. 235-245 (1978)

\bibitem{Breen-Messing (2001)} L. Breen and W. Messing,
``Differential Geometry of Gerbes'', preprint,
http://xxx.lanl.gov/abs/math.AG/0106083 (2001)

\bibitem{Brylinski (1993)} J.-L. Brylinski,
\underline{Loop Spaces, Characteristic Classes and Geometric Quantization},
Progress in Mathematics 107, Birkh\"auser, Boston (1993)

\bibitem{Carey-Mickelsson-Murray (1997)} A. Carey, J. Mickelsson and
M. Murray, 
``Index Theory, Gerbes, and Hamiltonian
Quantization'', Comm. Math. Phys. 183, p. 707-722 (1997)

\bibitem{Carey-Murray (1996)} A. Carey and M. Murray, 
``Faddeev's Anomaly and Bundle Gerbes'', Lett. Math. Phys. 37,
p. 29-36 (1996)

\bibitem{Dai (1991)} X. Dai, ``Adiabatic Limits, Nonmultiplicativity of 
Signature, and Leray Spectral Sequence'', J.
Amer. Math. Soc. 4, p. 265-321 (1991)

\bibitem{Dai-Zhang (1998)} X. Dai and W. Zhang, 
``Higher Spectral Flow'', J. of Funct. Anal. 157, p. 432-469 (1998)

\bibitem{Dold-Thom (1958)} A. Dold and R. Thom, 
``Quasifaserungen und Unendliche Symmetrische Produkte'', 
Ann. of Math. 67, p. 239--281 (1958)

\bibitem{Ekstrand-Mickelsson (2000)} C. Ekstrand and J. Mickelsson,
``Gravitational Anomalies, Gerbes, and Hamiltonian Quantization'', 
Comm. Math. Phys. 212, p. 613-624 (2000)

\bibitem{Faddeev (1984)} L. Faddeev, ``Operator Anomaly for the Gauss Law'',
Phys. Lett. 145B, p. 81-84 (1984)

\bibitem{Fegan-Gilkey (1985)} H. Fegan and P. Gilkey, 
``Invariants of the Heat Equation'', Pacific J. Math. 117, p. 233-254 (1985)

\bibitem{Hitchin (1999)} N. Hitchin,
``Lectures on Special Lagrangian Submanifolds'', 
in 
\underline{Winter School on Mirror Symmetry, Vector Bundles and 
Lagrangian Submanifolds},
AMS/IP Stud. Adv. Math. 23, 
Amer. Math. Soc., Providence, p. 151-182 (2001)

\bibitem{Lott (1994)} J. Lott, ``$\R/\Z$ Index Theory'', Comm. Anal.
Geom. 2, p. 279-311 (1994)

\bibitem{Melrose-Piazza (1997)} R. Melrose and P. Piazza,
``Families of Dirac Operators, Boundaries and the $b$-Calculus'',
J. Diff. Geom. 46, p. 99-180 (1997)

\bibitem{Miscenko-Fomenko (1979)}
A. Miscenko and A. Fomenko,
``The Index of Elliptic Operators over $C\sp{*}$-Algebras'',
Izv. Akad. Nauk SSSR Ser. Mat. 43, p. 831-859 (1979)

\bibitem{Pressley-Segal (1986)} A. Pressley and G. Segal,
\underline{Loop Groups}, Oxford University Press, Oxford (1986)

\bibitem{Quillen (1985)} D. Quillen, 
``Determinants of Cauchy-Riemann Operators on Riemann Surfaces'',
Funct. Anal. Appl. 19, p. 37-41 (1985)

\bibitem{Quillen (1985b)} D. Quillen, 
``Superconnections and the Chern Character'', Topology 24, p. 89-95 (1985)

\bibitem{Zumino-Wu-Zee (1984)} 
B. Zumino, Y. S. Wu and A. Zee, 
``Chiral Anomalies, Higher Dimensions, and Differential Geometry'', 
Nucl. Phys. B 239, p. 477-507 (1984)
\end{thebibliography}
\end{document}